\documentclass[12pt,leqno,a4paper]{amsart}
\usepackage{amssymb,enumerate}
\newcommand{\FF}{{\mathbb{F}}}
\newcommand{\QQ}{{\mathbb{Q}}}
\newcommand{\RR}{{\mathbb{R}}}

\newcommand{\bB}{{\mathbf{B}}}
\newcommand{\bC}{{\mathbf{C}}}
\newcommand{\bG}{{\mathbf{G}}}
\newcommand{\bH}{{\mathbf{H}}}
\newcommand{\bL}{{\mathbf{L}}}
\newcommand{\bN}{{\mathbf{N}}}
\newcommand{\bO}{{\mathbf{O}}}
\newcommand{\bP}{{\mathbf{P}}}
\newcommand{\bS}{{\mathbf{S}}}
\newcommand{\bT}{{\mathbf{T}}}
\newcommand{\bU}{{\mathbf{U}}}
\newcommand{\bZ}{{\mathbf{Z}}}

\newcommand{\fA}{{\mathfrak{A}}}
\newcommand{\fS}{{\mathfrak{S}}}

\newcommand{\Gal}{{\operatorname{Gal}}}

\newcommand{\Irr}{{\operatorname{Irr}}}
\newcommand{\Syl}{{\operatorname{Syl}}}
\newcommand{\GAP}{{\sf{GAP}}}

\newcommand{\GL}{{\operatorname{GL}}}

\newcommand{\SL}{{\operatorname{SL}}}
\newcommand{\PSL}{{\operatorname{PSL}}}
\newcommand{\GU}{{\operatorname{GU}}}
\newcommand{\SU}{{\operatorname{SU}}}

\newcommand{\PSU}{{\operatorname{PSU}}}
\newcommand{\Sp}{{\operatorname{Sp}}}

\newcommand{\SO}{{\operatorname{SO}}}

\newcommand{\tw}[1]{{}^{#1}\!}
\newcommand\hRu{{\widehat{R_u}}}
\newcommand\hlf{{\frac{1}{2}}}

\newcommand\bq{{\bar q}}
\newcommand\mn{{\!-\!}}
\newcommand\pl{{\!+\!}}
\newcommand\subn{{\,\triangleleft\triangleleft\,}}
\newcommand{\Sub}{{\operatorname{Sub}}}

\let\al=\alpha
\let\eps=\epsilon
\let\la=\lambda
\let\si=\sigma

\theoremstyle{theorem}
\newtheorem{thm}{Theorem}[section]
\newtheorem{lem}[thm]{Lemma}
\newtheorem{prop}[thm]{Proposition}
\newtheorem{cor}[thm]{Corollary}

\newtheorem{thmA}{Theorem}
\newtheorem{conjA}[thmA]{Conjecture}

\theoremstyle{definition}
\newtheorem{rem}[thm]{Remark}
\newtheorem{exmp}[thm]{Example}

\begin{document}

\title[Picky elements, subnormalisers, character correspondences]{Picky elements, subnormalisers, and\\ character correspondences}

\author{Gunter Malle}
\address{FB Mathematik, RPTU, Postfach 3049,
         67653 Kaisers\-lautern, Germany.}
\email{malle@mathematik.uni-kl.de}

\begin{abstract}
We gather evidence on a new local-global conjecture of Moret\'o and
Rizo on values of irreducible characters of finite groups. For this we study
subnormalisers and picky elements in finite groups of Lie type and determine
them in many cases, for unipotent elements as well as for semisimple elements
of prime power order. We also discuss subnormalisers of unipotent and
semisimple elements in connected and disconnected reductive linear algebraic
groups.
\end{abstract}

\thanks{The author gratefully acknowledges support by the DFG --Project-ID
286237555--TRR 195.}



\date{\today}

\maketitle


\section{Introduction}   \label{sec:intro}

In this paper we gather evidence for a new far-reaching conjecture of
Alex Moret\'o and Noelia Rizo on character values of finite groups.
Let $G$ be a finite group. For an element $x\in G$ let $\Irr^x(G)$ denote the
set of irreducible complex characters of $G$ that do not vanish at $x$.
Define the \emph{subnormaliser} of a subgroup $H$ of $G$ to be
$$S_G(H):=\{g\in G\mid H\subn\langle g,H\rangle \}$$
and for $x\in G$ let
$$\Sub_G(x):=\big\langle S_G(\langle x\rangle)\big\rangle.$$

The following conjecture was put forward by Moret\'o and Rizo \cite{MR24}:

\begin{conjA}[Moret\'o--Rizo]   \label{conj:AN}
 Let $G$ be a finite group and $p$ a prime. Then for any $p$-element $x\in G$
 there exists a bijection $f_x:\Irr^x(G)\to\Irr^x(\Sub_G(x))$ such that
 \begin{enumerate}
  \item[\rm(1)] $\chi(1)_p=f_x(\chi)(1)_p$, and
  \item[\rm(2)] $\QQ(\chi(x))=\QQ(f_x(\chi)(x))$.
 \end{enumerate}
\end{conjA}

Let us comment. (A much more thorough discussion of the conjecture and various
extensions of it is given in \cite{MR24}). First, observe that $\Sub_G(x)$
contains the normaliser of any Sylow $p$-subgroup $P$ of $G$ containing $x$.
In fact, as we show, it is generated by these. Thus, by the now proven McKay
conjecture there exist bijections between $\Irr_{p'}(\bN_G(P))$ and both
$\Irr_{p'}(G)$ and $\Irr_{p'}(\Sub_G(x))$, hence also between the latter two.
Since $p'$-degree characters do not vanish on any $p$-element, this would
form part of the required bijection~$f_x$. Now, in addition,
Conjecture~\ref{conj:AN} predicts a bijection on characters in 
$\Irr^x(G)\setminus\Irr_{p'}(G)$.

By Navarro's extension of the McKay
conjecture, there should even exist a $\Gal(\overline\QQ_p/\QQ_p)$-equivariant
bijection $\Irr_{p'}(G)\to\Irr_{p'}(\Sub_G(x))$, where $\QQ_p$ denotes the field
of $p$-adic numbers. Thus, and in view of the examples we discuss in this
paper, we are led to ask whether in the setting of Conjecture~\ref{conj:AN}
there exists a bijection $f_x$ such that moreover for any $\chi\in\Irr^x(G)$ we
have
 \begin{enumerate}
  \item[\rm(3)] $\chi(x)_p=f_x(\chi)(x)_p$, and
  \item[\rm(4)] $\QQ_p(\chi)=\QQ_p(f_x(\chi))$.
 \end{enumerate}
\medskip

Here for an algebraic number $\al\in\overline\QQ$ we write
$\al_p:=\big|$N$_{\QQ(\al)/\QQ}(\al)\big|_p^{1/[\QQ(\al):\QQ]}\in\RR_{\ge0}$
for its $p$-adic valuation.

Of course the conjecture is only meaningful when $\Sub_G(x)<G$.
An interesting special case is as follows: A $p$-element $x\in G$ is called
\emph{picky} if it belongs to a unique Sylow $p$-subgroup of $G$. Note that
$G$ has a picky $p$-element if and only if it does not have a \emph{redundant
Sylow $p$-subgroup} in the sense of \cite{MMM}. As an example, assume $G$ has
\emph{trivial intersection} (TI) Sylow $p$-subgroups. Then by definition any
non-identity $p$-element lies in a unique Sylow $p$-subgroup of $G$ and thus is
picky. This is the case, in particular, when $G$ has cyclic Sylow $p$-subgroups
of order~$p$. Now note that $x$ is picky in $G$ if and only if
$\Sub_G(x)=\bN_G(P)$, for $P\le G$ a Sylow $p$-subgroup of $G$ containing $x$
(see Corollary~\ref{cor:picky crit}).
\medskip

In this paper we undertake to test Conjecture~\ref{conj:AN} and Properties~(3)
and~(4) above in simple groups.
The paper is organised as follows. In Section~\ref{sec:basic} we collect some
basic results on picky elements and subnormalisers. In Section~\ref{sec:Lie
type unip} we classify the picky unipotent elements in groups of Lie type
and determine the subnormalisers of unipotent elements in most types. Based on
this we prove in Section~\ref{sec:conj type unip} the validity of
Conjecture~\ref{conj:AN} for unipotent elements of various families of groups
of Lie type. In Section~\ref{sec:Lie type ss} we classify picky semisimple
$p$-elements in groups of simply connected Lie type for all $p\ne2$. In the
final Section~\ref{sec:alg group} we discuss subnormalisers of unipotent and
semisimple elements in connected and disconnected reductive algebraic groups.
\medskip

\noindent
{\bf Acknowledgements:}
I thank Alex Moret\'o and Noelia Rizo for informing me of their conjectures and
for several enlightening conversations on the topic. I would also like to thank
Frank Himstedt for showing me the preprint \cite{Ya}, Martin van Beek for
providing the example in Remark~\ref{rem:vB}, Frank L\"ubeck for computing the
elements of order~9 in $E_7(2)$ and $E_8(2)$, and Mandi Schaeffer Fry for
comments on an earlier draft.

\section{Basic observations}   \label{sec:basic}
Let's make some easy observations about picky $p$-elements. Throughout $G$ is a
finite group and $p$ is a prime number.

\begin{lem}   \label{lem:centr}
 Let $P\le G$ be a Sylow $p$-subgroup and $x\in P$ picky in $G$. Then
 $\bN_G(\langle x\rangle)\le\bN_G(P)$.
\end{lem}

\begin{proof}
Let $g\in \bN_G(\langle x\rangle)$. Then $\langle x\rangle\le P^g$, so $x$ lies
in the Sylow $p$-subgroups $P$ and $P^g$. As $x$ is picky in $G$, we must have
$g\in \bN_G(P)$.
\end{proof}

\begin{lem}   \label{lem:abelian}
 Let $P\le G$ be a Sylow $p$-subgroup and $x\in P$. If both $P$ and $\bC_G(x)$
 are abelian then $x$ is picky.
\end{lem}

\begin{proof}
If $x\in P,P^g$ for some $g\in G$ then $P,P^g\le\bC_G(x)$, which has a unique
Sylow $p$-subgroup, being abelian, so $P=P^g$.
\end{proof}

\begin{lem}   \label{lem:normal}
 Let $H\le G$ with $[G:H]$ prime to $p$ and $x\in H$ a $p$-element.
 \begin{enumerate}[\rm(a)]
  \item If $x$ is picky in $G$ then $x$ is picky in $H$.
  \item If $H\unlhd G$ then $x$ is picky in $G$ if and only if it is picky
   in $H$.
 \end{enumerate}
\end{lem}

\begin{proof}
Part~(a) follows as any Sylow $p$-subgroup of $H$ is one of $G$, while in (b),
any $p$-element and any Sylow $p$-subgroup of $G$ is contained in~$H$.
\end{proof}

\begin{lem}   \label{lem:cover}
 Let $N\unlhd G$ where $N$ is either a $p$-group or central.
 Then a $p$-element $x\in G$ is picky if and only if $xN$ is picky in $G/N$.
\end{lem}

\begin{proof}
If $N$ is a $p$-group, the Sylow $p$-subgroups of $G/N$ are of the form $P/N$
for $P$ a Sylow $p$-subgroup of $G$, from which the assertion follows. If $N$
is central, then by the previous part we may assume it is a $p'$-group. Let
$x\in P_i$ for $P_i\in\Syl_p(G)$ with $P_1\ne P_2$. As $P_i$ is the unique
Sylow $p$-subgroup of $P_iN$ we also have $P_1N\ne P_2N$ and so $xN$ lies in
two distinct Sylow $p$-subgroups of $G/N$. The reverse direction is clear.
\end{proof}

Thus to study picky elements of a simple group $S$ we can instead consider a
covering group of $S$, or an extension of $S$ by a group of $p'$-automorphisms.
\medskip

Next, we collect some elementary properties of subnormalisers.

\begin{lem}   \label{lem:in Op}
 Let $H\le G$ and $x\in G$ a $p$-element with $\langle x\rangle\subn H$. Then
 $x\in\bO_p(H)$.
\end{lem}

\begin{proof}
Let $\langle x\rangle\unlhd N_1\unlhd\cdots\unlhd N_r=H$ be a subnormal series.
Clearly $x\in \bO_p(N_1)$, and since $\bO_p(N_i)$ is characteristic in $N_i$ we
have $\bO_p(N_i)\le\bO_p(N_{i+1})$ for all $i$, whence the claim.
\end{proof}

The following characterisation turns out to be very useful and could also be
taken as definition of subnormaliser for $p$-elements. It shows that the
subnormaliser of a $p$-element~$x$ is very closely related to the $p$-local
structure ``around $x$'':

\begin{prop}   \label{prop:gen NG(P)}
 Let $x\in G$ be a $p$-element. Then $\Sub_G(x)$ is generated by the
 normalisers of those Sylow $p$-subgroups of $G$ that contain $x$.
\end{prop}

\begin{proof}
Clearly $\langle x\rangle$ is subnormal in the normaliser of any Sylow
$p$-subgroup
containing~$x$, and hence $\Sub_G(x)$ contains all of these normalisers. For
the converse, first assume that $\langle x\rangle\subn H$ for some subgroup
$H\le G$. Then $x\in\bO_p(H)$ by Lemma~\ref{lem:in Op}, so $x$ is contained in
all Sylow $p$-subgroups of $H$. Let $P$ be a Sylow $p$-subgroup of $H$. By the
Frattini argument, we have $H=\bO^{p'}(H)\bN_H(P)$, hence $H$ is generated by
the normalisers in $H$ of its Sylow $p$-subgroups, each of which contains $x$.
\par
Now let $N\le G$ be any subgroup with a normal Sylow $p$-subgroup $Q$
containing~$x$. In particular, $x$ is subnormal in $N$. Let
$\langle x\rangle\unlhd Q_1\unlhd\ldots\unlhd Q\unlhd N$ be a subnormal series
in~$N$. If $Q$ is not a Sylow $p$-subgroup of $G$ there is a $p$-subgroup
$P\le G$ with $Q\unlhd P$ and $P>Q$. Since $Q\unlhd P$ then
$\langle x\rangle\unlhd Q_1\unlhd\ldots\unlhd Q\unlhd H$ is a subnormal series
in $H:=\langle N,P\rangle$. That is, any $N$ as above is contained in a
subgroup with a larger Sylow $p$-subgroup in which $\langle x\rangle$ is still
subnormal.
Hence, combining with the previous paragraph, any subgroup of $G$ in which $x$
is subnormal is contained in a subgroup generated by (subgroups of) Sylow
$p$-normalisers of $G$ that contain $x$. Our claim follows.
\end{proof}

See also Proposition~\ref{prop:unip alg grp} for an analogue for algebraic
groups. This shows (see also \cite[Thm~2.9]{MR24}):

\begin{cor}   \label{cor:picky crit}
 A $p$-element $x\in G$ is picky if and only if $\Sub_G(x)=\bN_G(P)$ for a
 Sylow $p$-subgroup $P$ of $G$ containing $x$.
\end{cor}

Recall that a $p$-subgroup $Q\le G$ is \emph{radical} if $Q=\bO_p(\bN_G(Q))$.
Since Sylow $p$-subgroups are clearly radical, we also obtain (see
\cite[Prop.~2.1(a)]{Ca90}):

\begin{cor}   \label{cor:radical}
 Let $x\in G$ be a $p$-element. Then $\Sub_G(x)$ is generated by the normalisers
 of the radical $p$-subgroups of $G$ that contain $x$.
\end{cor}

The subnormaliser of a $p$-element controls its fusion in a Sylow subgroup;
Moret\'o and Rizo had obtained a different proof of this fact based upon 
\cite[Prop.~2.1]{Ca90}.

\begin{lem}   \label{lem:fusion}
 Let $P$ be a Sylow $p$-subgroup of $G$ and $x\in P$. Then $\Sub_G(x)$ is
 generated by the elements $g\in G$ with $x^g\in P$.
\end{lem}

\begin{proof}
Let $x^g\in P$. Then $x\in\tw{g}P$ whence $\tw{g}P\le\Sub_G(x)$. Since all
Sylow $p$-subgroups of $\Sub_G(x)$ are conjugate, there is $h\in\Sub_G(x)$
with $\tw{g}P=P^h$ and so $hg\in\bN_G(P)\le\Sub_G(x)$, giving $g\in\Sub_G(x)$.
\par
Conversely, let $R:=\langle g\in G\mid x^g\in P\rangle$. Clearly,
$\bN_G(P)\le R$. Assume $x\in\tw{h}P$ for some $h\in G$. Then $x,x^h\in P$ and
thus $h\in R$. Thus, $R$ contains the normalisers of all Sylow $p$-subgroups of
$G$ containing $x$ and hence $\Sub_G(x)$ by Proposition~\ref{prop:gen NG(P)}.
\end{proof}

\begin{cor}   \label{cor:fusion}
 Let $x\in G$ be a $p$-element. Then $x$ is $G$-conjugate to $y\in\Sub_G(x)$
 if and only if $x,y$ are already $\Sub_G(x)$-conjugate. In fact, $\Sub_G(x)$
 is the smallest subgroup of $G$ containing both $\bC_G(x)$ and a Sylow
 $p$-subgroup of $G$ with this property.
\end{cor}

\begin{proof}
If $x,y\in\Sub_G(x)$ are $G$-conjugate, then up to conjugation in $\Sub_G(x)$
we may assume they lie in a common Sylow $p$-subgroup $P$ (of $\Sub_G(x)$).
Then the first claim follows by Lemma~\ref{lem:fusion}. For the second, let
$R\le G$ be a subgroup containing $\bC_G(x)$ and a Sylow $p$-subgroup $P$ of $G$
(which we may assume contains $x$) with the stated property. Let $g\in G$ with
$x^g\in P$ and thus $x,x^g$ are $R$-conjugate by assumption. That is, there
exists $r\in R$ with $x=x^{gr}$, whence $gr\in\bC_G(x)\le R$ and so $g\in R$.
Lemma~\ref{lem:fusion} then shows $\Sub_G(x)\le R$.
\end{proof}

\begin{rem}   \label{rem:vB}
It is tempting to ask whether the $p$-fusion system of the subnormaliser of
a $p$-element of $G$ is determined by the $p$-fusion system of $G$. But
this is not the case. I'm indebted to Martin van Beek for the following
counter-example: the 3-fusion categories of $M_{12}$ and $\PSL_3(3)$ are
known to be isomorphic. Yet, as we will show, $\PSL_3(3)$ has a (regular
unipotent) picky 3-element, while the corresponding 3-elements in $M_{12}$ have
centraliser isomorphic to $3\times\fA_4$ and subnormaliser $M_{12}$.
\end{rem} 

\begin{prop}   \label{prop:Burnside}
 Let $P\le G$ be a Sylow $p$-subgroup and $x\in P$. If $P$ is abelian then
 $\Sub_G(x)=\langle \bC_G(x),\bN_G(P)\rangle$.
\end{prop}

\begin{proof}
By Burnside's theorem, in this case $\bN_G(P)$ controls fusion of elements
in~$P$, so the claim follows from Corollary~\ref{cor:fusion}.
\par
Alternatively, let $H:=\langle \bC_G(x),\bN_G(P)\rangle$. Since $x$ is
subnormal in both
$\bC_G(x)$ and $\bN_G(P)$ we have $H\le \Sub_G(x)$. On the other hand, as $P$
is abelian, $\bC_G(x)$, and hence~$H$, contains all Sylow $p$-subgroups of $G$
containing~$x$. As $H$ contains $\bN_G(P)$, and all Sylow $p$-subgroups of $H$
are $H$-conjugate, it even contains the normalisers of all Sylow
$p$-subgroups of $G$ containing~$x$. Hence $\Sub_x(G)\le H$ by
Proposition~\ref{prop:gen NG(P)}, and we are done.
\end{proof}

Observe that Proposition~\ref{prop:gen NG(P)}, Lemma~\ref{lem:fusion},
Corollary~\ref{cor:fusion} and Proposition~\ref{prop:Burnside} hold more
generally with $x$ replaced by any $p$-subgroup
$H$ of~$G$ and $\bC_G(x)$ by $\bN_G(H)$, with identical proofs.

\begin{lem}   \label{lem:direct}
 Let $G=H_1\times\cdots\times H_r$ and $x=(x_1,\ldots,x_r)\in G$. Then
 $\Sub_G(x)=\Sub_{H_1}(x_1)\times\cdots\times\Sub_{H_r}(x_r)$. In particular,
 $x$ is picky in $G$ if and only if the projection of $x$ into each component
 $H_i$ is.
\end{lem}

\begin{proof}
Clearly, $g=(g_1,\ldots,g_r)\in S_G(x)$ if and only if $g_i\in S_{H_i}(x_i)$
for all $i$, showing the first claim. Since the Sylow $p$-subgroups of $G$ are
of the form $\prod_{i=1}^r P_i$, with Sylow $p$-subgroups $P_i\le H_i$, the
second assertion follows from the first using Corollary~\ref{cor:picky crit}.
\end{proof}

The next observation allows one to bound subnormalisers from above:

\begin{lem}   \label{lem:bound sub}
 Let $x\in G$ be a $p$-element with $x\in P\in\Syl_p(G)$. If $H\le G$
 contains $\bN_G(P)$, and $x\in H^g$ for $g\in G$ implies $H^g=H$, then
 $\Sub_G(x)\le H$. If moreover $\langle x\rangle\subn H$, then $\Sub_G(x)=H$.
\end{lem}

\begin{proof}
Let $g\in G$ be such that $x\in P^g$. Then $x\in P^g\le H^g$, so $H^g=H$ by
assumption, and moreover $\bN_G(P^g)=\bN_G(P)^g\le H^g=H$. Thus, $H$ contains the
normalisers of all Sylow $p$-subgroups containing $x$, whence the claim by
Proposition~\ref{prop:gen NG(P)}.
\end{proof}

The following result extends Lemma~\ref{lem:cover}:

\begin{lem}   \label{lem:mod unip}
 Let $N\unlhd G$, where $N$ is either a $p$-subgroup or central, and $x\in G$ a
 $p$-element. Then $N\le \Sub_G(x)$ and $\Sub_{G/N}(xN)= \Sub_G(x)/N$.
\end{lem}

\begin{proof}
If $N$ is central the claim is obvious. Assume $N$ is a $p$-group. Since
$\langle x,N\rangle$ is a $p$-group by assumption, we have
$\langle x\rangle\subn\langle x,N\rangle$ and so $N\le \Sub_G(x)$.
Suppose $gN\in S_{G/N}(xN)$, hence
$\langle xN\rangle\subn\langle gN,xN\rangle$. Then
$\langle x,N\rangle\subn\langle g,x,N\rangle$ and so
$\langle x\rangle\subn\langle g,x,N\rangle$, yielding $g\in S_G(x)$. The
reverse inclusion is analogous. See also \cite[Lemma~2.3]{Ca90}.
\end{proof}

Is is tempting to define $G$ has \emph{almost normal Sylow $p$-subgroups} if
$\Sub_G(x)=G$ for all $p$-elements $x\in G$. This is in some sense at the
opposite extreme from TI-Sylow $p$-subgroups, for which any $p$-element $x\ne1$
lies in exactly one Sylow $p$-subgroup, but the two notions coincide when the
Sylow $p$-subgroup is normal. An example with non-normal Sylow $p$-subgroup
for $p=2$ is the group {\tt SmallGroup(324,37)}, of the form $3^3.\fA_4$. In
fact, similar examples exist for all primes $p$. Non-solvable examples of
almost normal Sylow $p$-subgroups are given by the simple groups $\PSL_3(7)$,
$\PSL_3(13)$, $\PSU_3(5)$ and $\PSU_3(11)$ at $p=3$.

\section{Unipotent elements in groups of Lie type}   \label{sec:Lie type unip}
Here we investigate picky elements and subnormalisers for unipotent elements of
finite groups of Lie type. More precisely, we assume $\bG$ is a connected
reductive linear algebraic group over an algebraically closed field of
characteristic~$p$ and $F:\bG\to\bG$ is a Steinberg
map with (finite) group of fixed points $G:=\bG^F$ (see e.g.\ \cite{MT11}).
Observe that unipotent elements of $G$ are exactly the $p$-elements of~$G$. The
case of $p'$-elements will be discussed in Section~\ref{sec:Lie type ss}.

\subsection{Unipotent picky elements}
Let $T\le B$ be a maximal torus contained in a Borel subgroup of $G$. Thus $B$
is the group of $F$-fixed points of an $F$-stable Borel subgroup $\bB\le\bG$
with $F$-stable maximal torus $\bT$ where $T=\bT^F$. Then $B,N:=\bN_G(\bT)$
form a split BN-pair in~$G$ with Weyl group $W:=N/T$. Let $\Phi$
be the root system of this BN-pair, and $\Phi^+,\Delta\subset\Phi$ be the
positive respectively simple roots determined by $B$. Now $U:=\bO_p(B)$ is a
Sylow $p$-subgroup of $G$ and $B$ is its normaliser. Any element in $U$ can be
written uniquely as a product of root elements from the root subgroups $U_\al$,
for $\al\in\Phi^+$, contained in $U$ (see \cite[\S23]{MT11}).

\begin{prop}   \label{prop:BN}
 Let $\bG$ be connected reductive with Steinberg map $F$. Then a unipotent
 element $x\in G=\bG^F$ is picky if and only if, up to conjugation, $x$ is a
 product of root elements of $G$ in $B$ in which all simple roots do appear.
\end{prop}

\begin{proof}
Assume $x\in U=\bO_p(B)$ can be written as a product of root elements not
involving elements from $U_\al$ where $\al\in\Delta$. Let $s\in W$ be the
corresponding simple reflection and $\dot s\in \bN_G(\bT)$ a preimage. Then
$x^{\dot s}\in B$, that is, $x$ lies in the Borel subgroups $B$ and
$B^{\dot s}$ which are distinct, as $B^{\dot s}$ does not contain root elements
for the root $\al$. Hence $x$ lies in two different Sylow $p$-subgroups of $G$
and thus is not picky.

For the converse, assume $x$ lies in two Sylow $p$-subgroups, so (up to
conjugation) $x\in B$ and $x\in B^g$ for some $g\in G$. Writing $g$ in
Bruhat decomposition $g=u_1t\dot{w}u_2$ with $u_i\in U$, $t\in T$ and $w\in W$
we find $x\in (B\cap B^{w})^{u_2}$, so up to conjugation $x\in B\cap B^{w}$ for
some $1\ne w\in W$. Now by \cite[Prop.~2.5.9]{Ca} then $x\in U_{w_0w}$, with
$w_0\in W$ the longest element, and $U_{w_0w}$ is the product of the root
subgroups $U_\al$ for $\al\in\Phi^+\cap w(\Phi^+)$ by \cite[Prop.~2.5.16]{Ca}.
Since $w\ne 1$ there is some simple root $\al\in\Delta$ made negative by~$w$
and so the corresponding root element cannot occur in the unique expression
of $x$ as a product of root elements.
\end{proof}

We can now classify picky unipotent elements. Note that by
Lemmas~\ref{lem:normal} and~\ref{lem:cover} the precise isogeny type of our
connected reductive group $\bG$
does not matter since $[G:[G,G]]$ and $|\bZ(G)|$ are both prime to~$p$. Thus,
by virtue of Lemma~\ref{lem:direct} we are reduced to the case when $\bG$ is
simple, which we now assume.

\begin{thm}   \label{thm:unip picky}
 Let $\bG$ be simple simply connected with Steinberg map $F:\bG\to\bG$.
 A unipotent element $x\in G\setminus\{1\}$ is picky if and only if one of the
 following holds:
 \begin{enumerate}
  \item[\rm(1)] $x$ is regular unipotent;
  \item[\rm(2)] $G=\SU_{2n+1}(q)$ with $n\ge1$ and $x$ has Jordan block sizes $(2n,1)$;
  \item[\rm(3)] $G=\tw2B_2(2^{2f+1})$ with $f\ge0$ is a Suzuki group;
  \item[\rm(4)] $G={}^2G_2(3^{2f+1})$ with $f\ge0$ is a Ree group; or
  \item[\rm(5)] $G=\tw2F_4(2^{2f+1})$ with $f\ge0$ is a Ree group and $|\bC_G(x)|=2q^6$, for $q^2=2^{2f+1}$.
 \end{enumerate}
\end{thm}

\begin{proof}
By \cite[Prop.~5.1.3]{Ca}, a unipotent element $x\in\bG$ is regular if and only
if it lies in a unique Borel subgroup of $\bG$. Assume $x\in G$ lies in two
distinct Sylow $p$-subgroups, so in two Borel subgroups $B_1\ne B_2$ of $G$,
say $B_i=\bB_i^F$ for $F$-stable Borel subgroups $\bB_i$ of $\bG$, $i=1,2$.
Then $x$ lies in $\bB_1\ne\bB_2$ and hence is not regular unipotent. So regular
unipotent elements are picky.

For the converse, first assume $G$ is untwisted, so $F$ acts trivially on the
Weyl group
of $\bG$ and hence all root subgroups of $\bG$ in $\bB$ are $F$-stable. Let
$x\in U$ be not regular unipotent. Then, again by \cite[Prop.~5.1.3]{Ca}, it
can be written as a product of root elements in which at least one simple root
$\alpha$ does not occur, so it cannot be picky by Proposition~\ref{prop:BN}.

Now assume $F$ acts non-trivially on $W$ and hence on $\Delta$. Let $x\in U$
be picky. Then by Proposition~\ref{prop:BN} it is a product involving root
elements from all root subgroups for simple roots $\al\in\Delta$. Now by
\cite[Ex.~23.10]{MT11}, for example, the root subgroup $U_\al$ consists of
elements which are products of root elements of $\bG$ lying in an $F$-orbit of
simple roots of $\bG$, unless $\al$ is the image under the orbit map of a pair
of roots of $\bG$ forming a diagram of type $A_2,B_2$ or $G_2$. Thus, if we are
not in one of the latter cases, $x$ is a product of root elements of $\bG$
involving all simple roots and hence is regular.

It remains to discuss the cases for which the Dynkin diagram has a subgraph of
type $A_2,B_2$ or $G_2$ with non-trivial $F$-action. Thus $G$ is one of
$\SU_{2n+1}(q)$, $\tw2B_2(q^2)$, $^2G_2(q^2)$ or $\tw2F_4(q^2)$.
The rank~1 groups $\SU_3(q)$, $\tw2B_2(q^2)$ and $^2G_2(q^2)$ have TI Sylow
$p$-subgroups \cite[Prop.~2.3]{BLM}, so any of their non-trivial $p$-elements
is picky. For
$G=\tw2F_4(q^2)$ the list of unipotent class representatives in \cite{Sh75}
reveals that in addition to the regular classes, there are two
further cases, $u_{13},u_{14}$, involving root elements for both types of
simple roots, with centraliser as given in the statement, while all other
elements cannot be picky, by Proposition~\ref{prop:BN}. (Note that here the
elements from the two simple root subgroups are denoted $\al_1(t_1)\al_2(t_2)$
and
$\al_3(t)$ respectively.) Finally, for $G=\SU_{2n+1}(q)$ it remains to identify
picky elements $x$ that involve elements from all simple root subgroups of $G$,
but not from all simple root subgroups of $\bG$. Thus, expressed as a product
of root elements for $\bG$, $x$ involves root elements from $\bU_\al$ for all
$\al\in\{\al_1,\ldots,\al_{n-1},\al_n+\al_{n+1},\al_{n+2},\ldots,\al_{2n}\}$,
if the simple roots $\al_1,\ldots,\al_{2n}$ of $\bG$ are numbered along the
Dynkin diagram of $\bG$, and hence has Jordan block sizes $(2n,1)$ in the
natural representation.
\end{proof}

Note that the elements occurring in (2) and~(5) above are \emph{subregular},
that is, their centraliser dimension in $\bG$ is only~2 larger than that of
regular elements

\subsection{Subnormalisers of unipotent elements}
We keep the setting at the beginning of this section, so $G=\bG^F$ with
$\bG$ connected reductive and $F$ a Steinberg map.

Recall that the standard parabolic subgroups of $G$ are exactly the overgroups
of the Borel subgroup $B$, and are in natural bijection with the subsets
$\Gamma$ of $\Delta$ \cite[Prop.~12.2]{MT11}; we write $P_\Gamma$ for the
corresponding parabolic
subgroup. Let $x\in U$ be unipotent. Then the lattice of standard parabolic
subgroups $P$ of $G$ such that $x\in\bO_p(P)$ has a unique maximal element with
respect to inclusion, which we denote $P(x)$. Indeed, writing $x$ as a product
of root elements $u_\al\in U_\al$, with $\al\in\Phi^+$, we have $P(x)=P_\Gamma$
for $\Gamma:=\{\al\in\Delta\mid u_\al=1\}$. The Chevalley commutator formula
then shows that $P(x)=P(x^u)$ for all $u\in U$.

The following is a generalisation of the argument used to prove
Proposition~\ref{prop:BN}:

\begin{lem}
 Let $x\in U$ and $C=x^G$ the conjugacy class of $x$ in $G$. Assume that
 $P(x)\ge P(y)$ for all $y\in C\cap U$. Then $P=\Sub_G(x)$.
\end{lem}

\begin{proof}
Since $x$ is subnormal in $P(x)$ by definition, it suffices to see that any
$H\le G$ such that $\langle x\rangle\subn H$ lies in $P(x)$. Let $H$ be such a
subgroup. In particular $x\in\bO_p(H)$ by Lemma~\ref{lem:in Op}.
By the Borel--Tits theorem \cite[Thm~26.5]{MT11} there is a parabolic subgroup
$H'\ge H$ such that $x\in\bO_p(H')$, so we may assume $H$ is parabolic.
Furthermore, we may replace $H$ by any overgroup $H'$ with $x\in\bO_p(H')$. As
any such overgroup is again parabolic, we may assume $H$ is parabolic and
maximal with respect to $x\in\bO_p(H)$.   \par
Let $g\in G$ such that $H^g$ is a standard parabolic subgroup, say
$H^g=P_\Gamma$ for $\Gamma\subseteq\Delta$. Then
$x^g\in C\cap \bO_p(P_\Gamma)\subseteq C\cap U$ and hence $P_\Gamma\le P(x)$ by
assumption. In fact, by our choice of $H$ we have $P_\Gamma=P(x^g)$. Using the
Bruhat decomposition we can write $g=u\dot{w}b$ with $u\in U$,
$b\in B\le P_\Gamma$, $w\in W$, so $H^{uw}=P_\Gamma$ and
$x^{uw}\in\bO_p(P_\Gamma)$.
Since also $y:=x^u\in\bO_p(P_\Gamma)$ we have $y,y^w\in\bO_p(P_\Gamma)$.
Let $W_\Gamma\le W$ denote the Weyl group of the standard Levi subgroup of
$P_\Gamma$. Any $w\notin W_\Gamma$ sends at least one simple root in
$\Delta\setminus\Gamma$ into a negative root. On the other hand, since
$P_\Gamma=P(x^g)$ all simple roots from $\Delta\setminus\Gamma$ occur in $y$,
whence we can't have $\tw{w}y\in P_\Gamma$. Thus, $w\in W_\Gamma$, giving
$g\in P_\Gamma$ and so $H=\tw{g}P_\Gamma=P_\Gamma\le P(x)$. This proves our
claim.
\end{proof}

\begin{prop}
 Let $\bG$ be connected reductive with Steinberg map $F$. Let $x\in G=\bG^F$
 be unipotent. Then $\Sub_G(x)$ is a parabolic subgroup of $G$.
\end{prop}

\begin{proof}
By Proposition~\ref{prop:gen NG(P)}, $\Sub_G(x)$ contains a Sylow $p$-normaliser
of $G$, hence a Borel subgroup. Since all overgroups of $B$ are parabolic
subgroups, so is $\Sub_G(x)$.
\end{proof}

\begin{prop}   \label{prop:unip subnorm}
 Let $x\in B$ be unipotent. Suppose that $x$ can be written as a product of root
 elements in which no simple root in $\Delta_1\subseteq\Delta$ occurs. Then
 $P_{\Delta_1}\le\Sub_G(x)$. In particular, if $\Delta_1=\Delta$ then
 $\Sub_G(x)=G$.
\end{prop}

\begin{proof}
Let $\al\in\Delta_1$ be a simple root and $P_\al:=P_{\{\al\}}$ the
corresponding standard
parabolic subgroup of $G$. By assumption $x$ lies in the unipotent radical of
$P_\al$, so $\langle x\rangle$ is subnormal in $P_\al$. The claim follows since
$\langle P_\al\mid\al\in\Delta_1\rangle=P_{\Delta_1}$ (see
\cite[Prop.~12.2]{MT11}). The final assertion is now obvious.
\end{proof}

Subnormalisers in the simply laced, untwisted case are the easiest to determine:

\begin{prop}   \label{prop:GLn subnorm}
 Let $G$ be one of $\SL_n(q)$ ($n\ge2$), $\SO_{2n}^+(q)$ ($n\ge4$), $E_6(q)$,
 $E_7(q)$ or $E_8(q)$ and $x\in G$ be unipotent. Then $\Sub_G(x)=G$ if and only
 if $x$ is not regular.
\end{prop}

\begin{proof}
For regular elements this is Theorem~\ref{thm:unip picky}. For non-regular
elements, we first consider $G=\SL_n(q)$. For $n=2$ the only non-regular
unipotent
element is the identity, for $n=3$ the non-trivial non-regular unipotent
elements have Jordan type $(2,1)$ and thus are conjugate to a root element for
the non-simple positive root, so $\Sub_G(x)=G$ by Proposition~\ref{prop:unip subnorm}.
Now assume $n\ge4$. Let $x\in G$ be a non-regular unipotent element. Then $x$
has at least two Jordan blocks, so up to conjugation, $x$ can be written as a
product of simple root elements in which at least one simple root $\al$ is
missing. If this $\al$ lies at an end of the Dynkin diagram, then $x^{s_\al}$
is a conjugate in $B$ in which the simple root next to the end node is
missing. So we may assume $\al\ne\al_1,\al_{n-1}$. Let $P_1,P_2$ be the
end-node standard parabolic subgroups of $G$, corresponding to
$\Delta\setminus\{\al_1\}$, $\Delta\setminus\{\al_{n-1}\}$ respectively.
Let $L_i$ be the corresponding standard Levi subgroups, so $L_i\cong P_i/U_i$,
with $U_i=\bO_p(P_i)$ the unipotent radical of $P_i$. Then for $i\in\{1,2\}$,
the image $\bar x$
of $x$ in $L_i$ is a product of root elements not involving the simple root
$\al$ of $L_i$, hence not regular, so by induction $\Sub_{L_i}(\bar x)=L_i$.
(Note that $L_i\cong\GL_{n-1}(q)$, containing $\SL_{n-1}(q)$ as a normal
subgroup of $p'$-index.) As $U_i$ is a $p$-group, $\Sub_{P_i}(x)=P_i$
for $i=1,2$ by Lemma~\ref{lem:mod unip}. Since $\langle P_1,P_2\rangle=G$ by
\cite[Prop.~12.2]{MT11} this achieves the proof.
\par
Now assume $G$ is of one of the other listed types. Let $x\in G$ be non-regular
unipotent. By \cite[Prop.~5.1.3]{Ca}, up to conjugation, $x$ can be
written as a product of root elements in which at least one simple root
$\al$ is missing. Since the Dynkin diagram of $G$ has three end nodes, there
are at least two end node simple roots $\al_1,\al_2\in\Delta$ such that
$\al\in\Delta\setminus\{\al_i\}$, $i\in\{1,2\}$. Let $P_i$ be the corresponding
standard parabolic
subgroups of $G$. Since their Levi subgroups are of type $D_{n-1}$, $E_{n-1}$
or $A_3$, for which we know the result by induction, respectively by the first
part, we can argue exactly as before to see that $\Sub_{P_i}(x)=P_i$ for
$i=1,2$ and then $\Sub_G(x)\ge\langle P_1,P_2\rangle=G$.
\end{proof}

\subsection{Non-simply laced and twisted types}
In order to determine subnormalisers of unipotent elements in non-simply laced
type groups we first need to deal with the groups of rank~2. Note that for
the rank~1 groups $\SU_3(q)$, $\tw2B_2(q^2)$ and $^2G_2(q^2)$ the subnormalisers
of all non-trivial unipotent elements are Borel subgroups by
Theorem~\ref{thm:unip picky}. 

From now on, we number the simple roots in $\Delta$, for $\bG$ a simple group,
as in \cite[Tab.~9.1, Tab.~23.2]{MT11}, and we write
$P_i:=P_{\al_i}:=P_{\{\al_i\}}$ for $\al_i\in\Delta$.

\begin{prop}   \label{prop:rnk2 subnorm}
 Let $G$ be one of $B_2(q)$, $G_2(q)$ or $\tw3D_4(q)$ with $q=p^f$,
 or $\tw2F_4(q^2)$ with $q^2=2^{2f+1}$, and $x\in G$ unipotent. Then
 $\Sub_G(x)=G$ unless one of
 \begin{enumerate}
  \item[\rm(1)] $x$ is picky, where $\Sub_G(x)\sim_G B$;
  \item[\rm(2)] $G=B_2(q)$ with $p\ne2$ and $|\bC_G(x)|=2q^3(q+1)$, where
   $\Sub_G(x)\sim_G P_2$;
  \item[\rm(3)] $G=G_2(q)$ with $p\ne3$ and $|\bC_G(x)|=3q^4$, where
   $\Sub_G(x)\sim_G P_1$;
  \item[\rm(4)] $G=\tw3D_4(q)$ and $x$ is in class $D_4(a_1)$ with
   $|\bC_G(x)|=q^6$, where $\Sub_G(x)\sim_G P_2$; or
  \item[\rm(5)] $G=\tw2F_4(q^2)$ and $|\bC_G(x)|\in\{3q^{12},2q^8,4q^8\}$, where
   $\Sub_G(x)\sim_G P_1$.
 \end{enumerate}
\end{prop}

\begin{proof}
For all of these groups we know parametrisations of unipotent classes in
parabolic subgroups and expressions for class representatives in terms of
root elements. The picky elements have subnormaliser conjugate to $\bN_G(U)=B$,
by Corollary~\ref{cor:picky crit}. We consider the remaining classes.

Let us start with $G=B_2(q)$. Here representatives for the
unipotent conjugacy classes are given in \cite{En72,Sh82}. With the criterion
in Proposition~\ref{prop:unip subnorm} we see that for $p=2$ there cannot be a
class (apart from the regular ones) with $\Sub_G(x)<G$, and for $p>2$ only the
class with representative denoted $A_{22}$ in \cite{Sh82} could possibly have
that property. Representatives for the unipotent conjugacy classes of $G_2(q)$
are given in \cite{CR74,En76,EY86}, and again by Proposition~\ref{prop:unip
subnorm} it follows that only the class of $G_2(2^f)$ denoted $A_4$ in
\cite{EY86} and the class of $G_2(p^f)$, $p\ge5$, with representative $u_4$ in
\cite{CR74} could have $\Sub_G(x)<G$. For $G=\tw3D_4(q)$, consulting the tables
of class representatives in \cite{Ge91,Hi07}, we see that we only need to
consider the classes with representatives
$u_\beta(1)u_{2\alpha+\beta}(a)$ and $u_\al(1)u_{\alpha+\beta}(1)$. Now note
that conjugating $u_\beta(1)u_{2\alpha+\beta}(a)$ with the simple reflection
$s_\alpha$ gives an element not involving any simple root element. So we are
left with the class $D_4(a_1)$ with representative $u_\al(1)u_{\al+\beta}(1)$.
Finally, for $G=\tw2F_4(q^2)$ by \cite[Tab.~II]{Sh75} apart from the picky
cases identified in Theorem~\ref{thm:unip picky} only the representatives
$u_9,u_{10},u_{11},u_{12}$ involve one of the simple root elements.

It remains to determine the subnormaliser for the elements in (2)--(5).
By the tables in the cited literature, respectively in \cite{Hi04,Hi07} for
$\tw3D_4(q)$ and \cite{HH10} for $\tw2F_4(q^2)$, in all four groups we have the
following situation: exactly one of the unipotent radicals of the two maximal
standard parabolic subgroups contains an element $x$ as considered. Moreover,
the centraliser of $x$ in that parabolic subgroup is the same as in $G$. By
counting, in fact each such $x$ lies inside exactly one such maximal parabolic
subgroup, say~$P$. Since any parabolic subgroup contains a Borel subgroup,
hence a Sylow  normaliser, we conclude by Lemma~\ref{lem:bound sub} that
$\Sub_G(x)=P$.
\end{proof}

Observe that the unipotent elements in (2)--(4) above are again subregular.

\begin{lem}   \label{lem:Bn subnorm}
 Let $G=\SO_{2n+1}(q)$ ($n\ge3$) with $q$ odd and $x\in G$ be unipotent. Then
 $\Sub_G(x)\ne G$ if and only if either $x$ is regular, or if $x$ has Jordan
 form $(2n-1,1^2)$ and $|\bC_G(x)|=2q^{n+1}(q+1)$, in which case
 $\Sub_G(x)=P_n$.
\end{lem}

\begin{proof}
Let $P$ be the standard maximal parabolic subgroup of $G$ of type $A_{n-1}$,
that is, the stabiliser of a maximal isotropic subspace in the natural
representation, and $Q$ the standard maximal parabolic subgroup of type
$B_{n-1}$. We claim that if (up to conjugation) $\Sub_P(x)=P$ then
$\Sub_G(x)=G$. Indeed, assume $\Sub_P(x)=P$ and let $\bar x$
be the image of $x$ in the standard Levi subgroup $L\cong\GL_n(q)$ of $P$. By
Proposition~\ref{prop:GLn subnorm}, $\bar x$ is not regular, so (up to
conjugation) its image $y$ in the Levi factor of a standard parabolic subgroup
$L_1$ of $L$ of type $\GL_{n-1}(q)$ is not regular either. Now, $L_1$ is a
standard Levi subgroup of the standard Levi subgroup $M$ of $Q$ of type
$\SO_{2n-1}(q)$, and hence by induction $\Sub_M(y)=M$. (The induction base is
given by the case $2n-1=5$ from Proposition~\ref{prop:rnk2 subnorm}.) But then
$\Sub_Q(x)=Q$, and since $P$ is maximal in $G$ and $Q$ is distinct from $P$,
we have $\Sub_G(x)=G$.
\par
Now from the explicit class representatives given in \cite[\S4.1.2]{DLO} it can
be checked that unless $x$ is regular or has type
$V_\beta(2n-1)\oplus V_\beta(1)\oplus V(1)$ (in the notation of loc.\ cit.),
there is a conjugate $y\in P$ of $x$ whose image in $L$ is not regular, so
has $\Sub_P(y)=P$ and then $\Sub_G(y)=G$ by the first part. Finally assume $x$
has type $V_\beta(2n-1)\oplus V_\beta(1)\oplus V(1)$. Then $x$ has Jordan type
$(2n-1,1^2)$. In particular, writing $x\in U$ as a product of root elements,
all $u_\al$, $\al\in\Delta\setminus\{\al_n\}$, must be non-zero. But then the
only standard parabolic subgroups of $G$ that contain $x$ in their radical are
$B$ and $P_n$ (of type $B_1$).
The claim now follows with Lemma~\ref{lem:bound sub}.
\end{proof}

\begin{lem}   \label{lem:F4 subnorm}
 Let $G=F_4(q)$ and $x\in G$ be unipotent. Then $\Sub_G(x)=G$ unless $x$ is
 regular, or $q$ is odd and $x$ is in class $F_4(a_1)$, with $|\bC_G(x)|=2q^6$
 and representative $x_{24}$ in \cite[Tab.~5 resp.~~6]{Sho74}, where
 $\Sub_G(x_{24})=P_{\{3,4\}}$.
\end{lem}

\begin{proof}
Representatives for the unipotent conjugacy classes of $G$ as well as for a
Levi subgroup $L$ of $G$ of type $B_3$ were determined in \cite{Sh74,Sho74}.
First assume $p=2$. We claim that only the regular unipotent elements in $L$
have subnormaliser strictly smaller than~$L$. Indeed, seven of the class
representatives given in \cite[Prop.~2.1]{Sh74} do not involve any simple root
element, and the other three non-regular ones can be conjugated by a simple
reflection to elements $x$ which only involve the 3rd simple root of $L$.
Then $x$ has trivial image in the standard Levi subgroup of $L$ of type $B_2$,
and image involving just one simple root in the standard Levi subgroup of
type~$A_2$, so $\Sub_L(x)=L$ by Proposition~\ref{prop:unip subnorm}.

Now we discuss the unipotent class representatives of $G$ from
\cite[Thm~2.1]{Sh74}. The ones labelled $x_0$ through $x_{19}$ do not involve
any simple root element, and similarly for $x_{27},x_{28}$. For elements
$x\in\{x_{22},\ldots,x_{26}\}$ only one of the short simple roots occurs, so
the images of $x$ in standard Levi subgroups of type $A_2$ and $B_3$ are
non-regular and we conclude by Proposition~\ref{prop:unip subnorm} and the
previous paragraph. Similarly for $x=x_{20},x_{21}$ the image in a Levi of
type~$C_3$ is trivial, and in a Levi of type $B_3$ is non-regular. By
Lemma~\ref{lem:Bn subnorm} this leaves
us with the elements $x_{29}$ and $x_{30}$, of type $F_4(a_1)$. The given
representatives show that these elements are regular in the subsystem
subgroup of type $B_4$. Conjugating $x=x_{29}$ or $x_{30}$ with the simple
reflection not belonging to the $B_3$-subsystem, we obtain an element which
only involves simple root elements for the 2nd and 3rd simple root. By the first
paragraph, the subnormaliser of that element contains the standard parabolic
subgroup of type~$B_3$, but by Proposition~\ref{prop:GLn subnorm} also the one
of type $A_2$, hence it is all of $G$.

If $p$ is odd, Proposition~\ref{prop:rnk2 subnorm} shows that apart from the
regular classes, only the unipotent class of $L$ with representative $z_8$
(in the notation of \cite[Tab.~3]{Sho74}) has proper subnormaliser. From
\cite[Tab.~5 and~6]{Sho74} we conclude that again at most the classes of~$G$
of type $F_4(a_1)$ with representatives $x_{23},x_{24}$ might have a
proper subnormaliser. Conjugating $x=x_{23},x_{24}$ by the simple reflection
not in the $B_3$-subsystem we obtain elements whose image in $L$ equals
$z_7,z_8$ respectively. As observed above, $\Sub_L(z_7)=L$ and so
$\Sub_G(x_{23})=G$ by the usual argument. On the other hand, we have
$\Sub_L(z_8)$ is the standard parabolic subgroup corresponding to the 3rd node
of the diagram of $B_3$, which is also the 3rd node of the diagram of $F_4$.
Thus the set of standard parabolic subgroups of $G$ containing $x_{24}$ in their
radical has the unique maximal element $P_{\{3,4\}}$, and so
$\Sub_G(x_{24})=P_{\{3,4\}}$ by Lemma~\ref{lem:bound sub}.
\end{proof}

\begin{lem}   \label{lem:Su even subnorm}
 Let $G=\SU_{2m}(q)$ and $x\in G$ be unipotent. Then $\Sub_G(x)=G$ unless $x$
 is regular, or $x$ has Jordan type $(2m-1,1)$, when $\Sub_G(x)\sim_GP_m$.
\end{lem}

\begin{proof}
There is nothing to show for $m=1$, so assume $m\ge2$. Let $x\in G$ be
unipotent and denote by $\la_1\ge\la_2\ge\cdots$ the lengths of its Jordan
blocks. According to the normal forms given in \cite[4.1.3]{DLO}, $x$ has a
conjugate in $B$ in which the entry at position $(i,i+1)$ is zero for at least
$2m-\la_1$ indices $1\le i \le m$. Assume $2m-\la_1>1$. Then the image of~$x$
in the standard Levi subgroup $\GL_m(q^2)$ is not regular and thus $\Sub_G(x)$
contains the corresponding maximal parabolic subgroup by
Propositions~\ref{prop:unip subnorm} and~\ref{prop:GLn subnorm}. Furthermore,
the image of ~$x$ in
the standard Levi subgroup $\GU_{2m-2}(q)$ also has at least 2 zero entries
directly above the main diagonal, unless $\la_1=\la_2$. In the latter case
we must have $2m=4$, and then the image of $x$ in $\GU_{2m-2}(q)=\GU_2(q)$ is
trivial. So in either case, by induction $\Sub_G(x)$ also contains this end
node parabolic subgroup and thus equals $G$.

This only leaves the cases when
$2m-\la_1\le1$, that is, $x$ has Jordan type $(2m)$ or $(2m-1,1)$. The first
of these is the class of regular unipotent elements. If $x\in U$ has a Jordan
block of length~$2m-1$, then it must involve elements from all simple root
subgroups except possibly for the 'middle' one. In particular, the only
standard parabolic subgroups whose unipotent radical could contain $x$ are $B$
and $P_m$. It is easy to see that $\bO_p(P_m)$ contains elements with a Jordan
block of length $2m-1$; then Lemma~\ref{lem:bound sub} shows $\Sub_G(x)=P_m$.
\end{proof}

The subnormalisers of unipotent elements in $\Sp_{2n}(q)$ ($n\ge3$),
$\SU_{2m+1}(q)$ ($m\ge3$), $\SO_{2n}^-(q)$ ($n\ge4$) and $\tw2E_6(q)$ seem more
involved and we will not discuss them here. For example, in $\Sp_{2n}(q)$
the number of unipotent classes with proper subnormaliser seems to increase
with the rank~$n$.

\section{On Conjecture~\ref{conj:AN} for the defining prime}   \label{sec:conj type unip}
We now use our results on subnormalisers to verify Conjecture~\ref{conj:AN} for
unipotent elements of groups of Lie type of rank at most~2. We keep the
notation and setting from the beginning of Section~\ref{sec:Lie type unip}.
Of course only the elements $x\in G$ with $\Sub_G(x)<G$ are of interest. We
write \emph{Conjecture~\ref{conj:AN}$^+$} to include Properties~(3) and~(4)
from the introduction. Our first result is for groups of arbitrary rank.

\begin{prop}   \label{prop:good}
 Assume that $p$ is a good prime for $\bG$ and $\bZ(\bG)$ is connected. Then
 Conjecture~\ref{conj:AN}$^+$ holds for regular unipotent elements of $G$.
\end{prop}

\begin{proof}
Let $x\in B\le G$ be regular unipotent. Under our assumptions on $\bG$ and $p$,
by the theorem of Green--Lehrer--Lusztig \cite[Cor.~8.3.6]{Ca} the irreducible
characters of $G$ that do not vanish on~$x$ are exactly those of degree prime
to~$p$, there are $|\bZ(G)|q^l$ of these, where $l$ denotes the semisimple rank
of $\bG$, and they all take value~$\pm1$ on $x$. This implies also that $x$
must be rational. Now $\Sub_G(x)=B$ as $x$ is picky by
Theorem~\ref{thm:unip picky}. By the proven McKay conjecture for groups of Lie
type in defining characteristic \cite{Mas} there is the same number of
irreducible characters of $B$ of $p'$-degree. Moreover, as $x$ is picky, by
Lemma~\ref{lem:centr} it must also be rational in $B$ and
$|\bC_B(x)|=|\bC_G(x)|$. This implies that the $|\bZ(G)|q^l$
irreducible characters of $B$ not vanishing on $x$ must also take values $\pm1$.
\end{proof}

This of course leaves open the case of non-connected centre, as well as that of
bad primes; we'll discuss a few examples for the latter in which complete
character tables are known.

\begin{prop}   \label{prop:G2}
 Conjecture~\ref{conj:AN}$^+$ holds for $G_2(q)$, $q=p^f$, at the prime $p$.
\end{prop}

\begin{proof}
Let $x\in G:=G_2(q)$ be a $p$-element, hence unipotent. There is nothing to
prove when $\Sub_G(x)=G$, so we are in one of the cases of
Proposition~\ref{prop:rnk2 subnorm}. If $x$ is regular and $p$ is good for $G$
the claim follows by Proposition~\ref{prop:good}. So we need to consider
regular elements for the bad primes $p=2,3$, and the subregular unipotent class
for any $p\ne3$.

The character tables for $G=G_2(2^f)$ and a Borel subgroup $B$ of $G$ were
determined by Enomoto--Yamada \cite{EY86}. First let $x\in B$ be regular (there
are
two such classes) so $\Sub_G(x)=B$. According to Tables~I and~IV of that paper
both $\Irr^x(G)$ and $\Irr^x(B)$ consist of $q^2$ characters of $p'$-degree and
4 characters with $p$-part of the degree equal to~$q/2$. The values of the
former on $x$ are $\pm1$, and $\pm q/2$ for the latter in both groups, and their
rationality properties agree (for the $p'$-characters this follows by
\cite{Ru21}), so Conjecture~\ref{conj:AN}$^+$ holds for $x$.

The character table of $G=G_2(3^f)$ and of its Borel subgroup can be found in
\cite[Tab.~I and~VII]{En76}. The group $G$ contains three classes of regular
unipotent elements~$x$, where again $\Sub_G(x)\sim_GB$. Here, both $\Irr^x(G)$
and $\Irr^x(B)$ consist of $q^2$ characters of $p'$-degree and 6 characters with
$p$-part of the degree equal to $q/3$. The values of the $p'$-characters on~$x$
are all $\pm1$. The other characters take values $\pm 2q/3$ (4 times) and
$\pm q/3$ (twice) on one of the classes, on the other two the values are
$\pm q/3$ (4 times) and $q/3+\zeta_3^iq$ for $i=1,2$, where $\zeta$ is a
primitive third root of unity, both for $G$ and for $B$. Again,
Conjecture~\ref{conj:AN}$^+$ is seen to hold.

Now let $x\in G=G_2(2^f)$ be in the subregular unipotent class from
Proposition~\ref{prop:rnk2 subnorm}(3) where $\Sub_G(x)\sim_{P_1}P_1$, the
first maximal
standard parabolic subgroup. In Tables~\ref{tab:G2} and~\ref{tab:G2 P} we have
extracted from \cite{EY86} the values of all irreducible characters of $G$ and
$P_1$  (denoted $P_a$ in \cite{EY86}) \emph{of even degree} that do not vanish
on $x$. (The class representative is denoted $A_4$ in $G$ and $A_5$ in $P_1$ in
\cite{EY86}.) Here
$\eps\in\{\pm1\}$ is such that $q\equiv\eps\pmod3$. Visibly, there is a map as
required in Conjecture~\ref{conj:AN} which even preserves values up to sign.
For the characters of odd degree, such a map additionally preserving character
fields over $\QQ_p$ exists by \cite{Ru21}, so it satisfies
Conjecture~\ref{conj:AN}$^+$.
For $p\ge5$ the character table of $G_2(p^f)$ was found by Chang and Ree
\cite{CR74}, the one for $P_2$ by Yamada \cite{Ya}, and the same sets of values
arise.
\end{proof} 

\begin{table}[htbp]
\caption{Character values for $G_2(2^f)$ \ldots}   \label{tab:G2}
$\begin{array}{c|cccccccc}
 & \theta_1,\theta_1'& \theta_2,\theta_2'& \theta_3& \theta_4,\theta_9(\pm1)& \theta_8& \chi_2(k)& \chi_2'(k)\\
 \#& 2& 2& 1& 3& 1& \hlf(q\mn3\mn\eps)& \hlf(q\mn1\pl\eps)\\
\hline
 A_4& \frac{1}{6}q(\eps q-1)& -\hlf q(\eps q-1)& \frac{1}{3}q(\eps q-1)& \frac{1}{3}q(\eps q+2)& q& q& -q\\
\end{array}$
\bigskip

\caption{\ldots and for $P_1$}   \label{tab:G2 P}
$\begin{array}{c|ccccc}
 & \theta_3(\pm1)& \theta_2(\pm1)& \theta_4& \theta_5,\theta_6(\pm1)& \chi_3(k)\\
 \#& 2& 2& 1& 3& q\mn1\\
\hline
 A_5& \frac{1}{6}q(\eps q-1)& -\hlf q(\eps q-1)& \frac{1}{3}q(\eps q-1)& \frac{1}{3}q(\eps q+2)& q\\
\end{array}$
\end{table}

\begin{prop}   \label{prop:3D4}
 Conjecture~\ref{conj:AN}$^+$ holds for $\tw3D_4(q)$, $q=p^f$, at the prime $p$.
\end{prop}

\begin{proof}
Let $x\in G:=\tw3D_4(q)$ be unipotent. By Proposition~\ref{prop:rnk2 subnorm}
we know $\Sub_G(x)=G$ unless $x$ is regular or subregular. For regular elements
the claim follows by Proposition~\ref{prop:good} for odd $q$ as $G$ can be
constructed as the $F$-fixed points of a simple group $\bG$ of adjoint
type $D_4$ whose only bad prime is $p=2$.

For $p=2$, the character table of $G$ was computed in \cite{DM87} and the
one of a Borel subgroup in \cite[Tab.~A.6]{Hi07}. For both groups there are
$q^4$ characters of odd degree, all taking value $\pm1$ on regular unipotent
elements, and four characters of degree divisible by~$q^3$, taking values
$\pm q^2/2$.

Now we consider the elements $x$ in the subregular unipotent class $D_4(a_1)$
of $G$ with $\Sub_G(x)=P_2$, the second maximal standard parabolic subgroup.
Here, by the tables in \cite{DM87,Hi04,Hi07} both $G$ and $P_2$ possess $q^3$
irreducible characters of degree divisible by~$p$ (in fact, precisely by $q$)
not vanishing on $x$, and all of them take value $\pm q$ on $x$.
Conjecture~\ref{conj:AN}$^+$ now follows, using \cite{Ru21} for the character
fields of $p'$-characters.
\end{proof} 

The Suzuki groups $\tw2B_2(2^{2f+1})$ possess TI Sylow 2-subgroups. Moret\'o and
Rizo checked their conjectures for the 2-elements of this group \cite{MR24}.
A further interesting case is given by the Ree groups ${}^2G_2(3^{2f+1})$ whose
Sylow 3-subgroups are also TI, so all non-identity 3-elements are picky.

\begin{prop}   \label{prop:2G2}
 Conjecture~\ref{conj:AN}$^+$ holds for ${}^2G_2(q^2)$, $q^2=3^{2f+1}$, at
 $p=3$.
\end{prop}

\begin{proof}
The character table of $G:={}^2G_2(q^2)$ was found by Ward \cite{Wa66} while
that of a Sylow 3-normaliser was computed by van der Waall \cite[p.~173]{vdW}. 
In Table~\ref{tab:2G2} we reproduce the values of those characters in $\Irr(G)$
respectively $\Irr(B)$ of degree divisible by $p$ not vanishing on some
non-trivial unipotent element.

\begin{table}[htb]
\caption{Some character values for ${}^2G_2(q^2)$ and for $B$ on unipotent classes.}   \label{tab:2G2}
$\begin{array}{c|cccc||c|cccc}
  {}^2G_2(q^2)& 1& X& T& YT^i& \qquad B\qquad&  1& X& T& YT^i\\
 \hline
 \overline{\xi_5},\xi_7&  \hlf \bq\Phi_1\Phi_6''& -\hlf(q^2+\bq)&
  \bq b& \bq\zeta_3^i& 
  \psi_3,\overline{\psi_5}& \hlf \bq\Phi_1& \hlf \bq\Phi_1& \bq b& \bq\zeta_3^i \\
\overline{ \xi_6},\xi_8&  \hlf \bq\Phi_1\Phi_6'& \hlf(q^2-\bq)&
 \bq b& \bq\zeta_3^i&
  \psi_4,\overline{\psi_6}& \hlf \bq\Phi_1& \hlf \bq\Phi_1& \bq b& \bq\zeta_3^i \\
 \overline{\xi_9},\xi_{10}&  \bq(q^4-1)& -\bq& 2\bq b& -\bq\zeta_3^i&
  \psi_1,\overline{\psi_2}& \bq\Phi_1& \bq\Phi_1& 2\bq b& -\bq\zeta_3^i\\
 \xi_4&  q^2\Phi_6& q^2& .& .&   \psi&  q^2\Phi_1& -q^2& .& .\\
\end{array}$
\end{table}

Here $\bq:=q/\sqrt{3}$, $\Phi_1:=q^2-1$, $\Phi_6':=q^2-3\bq+1$,
  $\Phi_6'':=q^2+3\bq+1$, $\Phi_6:=\Phi_6'\Phi_6''=q^4-q^2+1$,
  $b:=(-1+\sqrt{-3}\bq)/2$, $\zeta_3:=(-1+\sqrt{-3})/2$ and $i\in\{1,2\}$.

A quick check shows that Conjecture~\ref{conj:AN} holds for all (picky) classes,
and furthermore a bijection preserving character fields over $\QQ_p$ exists
using \cite{J22}.
\end{proof}

The next example is quite interesting as there exists an element $x$ with
$\Sub_G(x)<G$ for which $\Irr^x(G)$ contains characters of six different heights
and yet all $p$-parts are preserved under a suitable bijection.

\begin{prop}   \label{prop:2F4}
 Conjecture~\ref{conj:AN}$^+$ holds for $\tw2F_4(q^2)$, $q^2=2^{2f+1}$, at
 $p=2$.
\end{prop}

\begin{proof}
We need to consider the regular classes, the picky subregular classes from
Theorem~\ref{thm:unip picky}(5) as well as the classes with proper
subnormaliser in Proposition~\ref{prop:rnk2 subnorm}(5). The values of
unipotent characters of $G$ are given in \cite{Ma90}, the values of all other
characters at least on unipotent elements can be found in \cite{Chv}. The
character table of a Borel subgroup of $G$ was computed in \cite{HH09}.
For $x$ in one of the four regular unipotent classes or the two
subregular picky classes, $\Irr^x(G)$ and $\Irr^x(B)$ both consists of $q^4$
characters of odd degree, $2q^2$ characters whose degree has 2-part
$\sqrt2q/2$, and 8 characters whose degree has 2-part~$q^4/4$. These take values
$\pm1$, $q\sqrt{-2}/2$, and $\pm q^2/2$, $\pm \sqrt{-1}q^2/2$ on regular
elements, and values $\pm1$, $q\sqrt{2}/2$ respectively $\pm\sqrt{2}q^3/4$ on
the subregular elements, both for $G$ and~$B$. A bijection preserving character
fields over $\QQ_p$ was given in \cite{J22}.

The parabolic subgroup $P_1$ also possesses $q^4$ characters of odd degree by
\cite{HH10}. For elements $x$ which are either picky or lie in the classes in
Proposition~\ref{prop:rnk2 subnorm}(5), the $p$-parts of the values of
characters in $\Irr^x(G)$ and $\Irr^x(P_1)$ of even degree are given in
Tables~\ref{tab:2F4} and~\ref{tab:2F4 NGP}, ordered by increasing $p$-part.
Here, we have set $\bq:=q/\sqrt{2}$ and $i:=\sqrt{-1}$. The 2-part of
$\bq^4\pm2i\bq^3$ occurring in either table is actually $2\bq^3$ for $q^2\ne8$
and $2\sqrt{2}\bq^3$ for $q^2=8$.

\begin{table}[htb]
\caption{$2$-Parts of character values for $\tw2F_4(q^2)$ \ldots}   \label{tab:2F4}
$\begin{array}{c|c|ccccccccc}
  & \#& 1& u_9& u_{10}& u_{11}=u_{12}^{-1}& u_{13-14}& u_{15-18}\\
 \hline
 \chi_{2,3,23,24}& 2q^2& \bq& \bq& \bq& \bq& \bq& \bq\\
 \chi_{4,27,30,33}& q^2& q^2& q^2& q^2& q^2& .& .\\
 \chi_{5,6,8,9}& 4& \bq^4& \bq^4& \bq^4& \bq^4& \bq^3& \bq^2\\
 \chi_{11-14}& 4& \bq^4& \bq^4& \bq^4& (\bq^4\pm2i\bq^3)_2& \bq^3& \bq^2\\
 \chi_{7,10}& 2& 2\bq^4& 2\bq^4& 2\bq^4& 2\bq^4& .& .\\
 \chi_{15-17}& 3& q^4& 2q^4& .& .& .& .\\
\end{array}$
\bigskip

\caption{\ldots and for $P_1$}   \label{tab:2F4 NGP}
$\begin{array}{c|c|ccccccccc}
  & \#& 1& c_{1,33}& c_{1,34}& c_{1,35}=c_{1,36}^{-1}& c_{1,37-38}& c_{1,39-42}\\
 \hline
 \chi_{7-8}(k),\chi_{9-10}& 2q^2& \bq& \bq& \bq& \bq& \bq& \bq\\
 \chi_2(k),\chi_{11}& q^2& q^2& q^2& q^2& q^2& .& .\\
 \chi_{21-24}& 4& \bq^4& \bq^4& \bq^4& \bq^4& \bq^3& \bq^2\\
 \chi_{14-17}& 4& \bq^4& \bq^4& \bq^4& (\bq^4\pm2i\bq^3)_2& \bq^3& \bq^2\\
 \chi_{13,25}& 2& 2\bq^4& 2\bq^4& 2\bq^4& 2\bq^4& .& .\\
 \chi_{18-20}& 3& q^4& 2q^4& .& .& .& .\\
\end{array}$
\end{table}

The tables show that there is a bijection as in Conjecture~\ref{conj:AN}, and it
can be checked it moreover preserves character fields over $\QQ_p$.
\end{proof}

\begin{prop}   \label{prop:Sp4}
 Conjecture~\ref{conj:AN}$^+$ holds for $\Sp_4(q)$, $q=p^f$, at the prime $p$.
\end{prop}

\begin{proof}
Here, Proposition~\ref{prop:good} does not apply as for odd $p$ the centre of
the algebraic group $\bG=\Sp_4$ is disconnected, and the prime~$p=2$ is bad for
$\bG$. For $p=2$ the character tables of $G$ and $B$  were calculated by Enomoto
\cite[Tab.~I and~IV]{En72}. There are $q^2$ irreducible characters of odd
degree, taking value $\pm1$ on a regular unipotent element $x$, and 4 characters
of degree divisible by $q/2$ taking values $\pm q/2$, both for $G$ and for $B$.

The character table of $G=\Sp_4(q)$ for odd $q$ was determined by Srinivasan
\cite{Sr68} and for~$B$ by Yamada \cite[Tab.~I-2]{Ya}. The characters not
vanishing on a regular unipotent element $x\in G$ for both $G$ and $B$ are all
of $p'$-degree, and $q(q-1)$ of them take values $\pm1$ on $x$, the other $4q$
take values $\pm(1\pm\sqrt{q^*})/2$, where $q^*:=(-1)^{(q-1)/2}q$.

Finally, for $p$ odd and $x$ the subregular element as in
Proposition~\ref{prop:rnk2 subnorm}(2) we have $\Sub_G(x)=P_1$. The character
table of $P_1$ can be found
in \cite[II-2]{Ya}, where the corresponding class is denoted $A_3(1)$. The
characters of degree divisible by $p$ take values $\pm q$ ($q+3$ times) and
$\pm2q$ ($(q-1)/2$ times) on $x$, for both $G$ and $P_1$.
\end{proof}

\section{Picky semisimple elements in groups of Lie type}   \label{sec:Lie type ss}
We now turn to the investigation of semisimple elements of groups of Lie type,
that is to say, elements of (prime power) order different from the defining
characteristic. We stay in the setting from the beginning of
Section~\ref{sec:Lie type unip}, with $\bG$ connected reductive and $G=\bG^F$ a
finite groups of Lie type, and let $\ell\ne p$ be a prime. To investigate
the existence of picky $\ell$-elements in simple groups $S$, by
Lemma~\ref{lem:cover}, we may consider $\bG$ of simply connected type such that
$G/\bZ(G)\cong S$.

\subsection{The case of abelian Sylow $\ell$-subgroups}
We first assume that $F$ is a Frobenius map with respect to an
$\FF_q$-structure, so $G$ is not a Suzuki or Ree group. Recall then (see e.g.\ 
\cite[\S3.5]{GM20}) that for any integer $d\ge1$, $(\bG,F)$ has Sylow $d$-tori,
that is, $F$-stable tori $\bS$ of $\bG$ whose order polynomial is the maximal
$\Phi_d$-power dividing the order polynomial of $(\bG,F)$. Then
$W_d:=\bN_G(\bS)/\bC_G(\bS)$ is the \emph{relative Weyl group} of the Sylow
$d$-torus. We write $e_\ell(q)$ for the order of $q$ modulo~$\ell$, respectively
modulo~4 if $\ell=2$.


\begin{lem}   \label{lem:cases}
 In the above setting, if $\ell$ divides $\Phi_d(q)$ for a unique cyclotomic
 factor $\Phi_d$ occurring in the order polynomial of $(\bG,F)$, then Sylow
 $\ell$-subgroups of $G$ are abelian. If $\bG$ is simple, the converse holds.
\end{lem}

This is shown in \cite[Prop.~2.2]{Ma14}. Note that we do not claim that the
simple groups $S=G/\bZ(G)$ satisfy this dichotomy; counter-examples occur for
$\SL_2(q)$ at $\ell=2$ and $\SL_3(\pm q)$ at $\ell=3$. The following is an
analogue of Proposition~\ref{prop:gen NG(P)}:

\begin{prop}   \label{prop:gen norm syl}
 In the above setting, assume that $\ell$ divides a unique cyclotomic factor
 $\Phi_d(q)$ in the order polynomial of $(\bG,F)$. Let $x\in G$ be an
 $\ell$-element. Then $\Sub_G(x)$ is generated by the $F$-fixed points of the
 normalisers of the Sylow $d$-tori of $\bG$ containing $x$.
\end{prop}

\begin{proof}
Let $x\in G$ be an $\ell$-element, and $P$ a Sylow $\ell$-subgroup of $G$
containing $x$. Then $P$ is abelian by Lemma~\ref{lem:cases} and lies in
a Sylow $d$-torus $\bS$ of $\bG$. Conversely, if $\bS$ is a Sylow $d$-torus
containing~$P$, then $\bS\le\bC_\bG(P)$, hence even
$\bS\le\bH:=\bC_\bG^\circ(P)$. By inspection of the order formulae
\cite[Tab.~24.1]{MT11} our condition on $\ell$ implies that $\ell$ does not
divide $|W|$. Thus, $\ell$ is not a torsion prime for $\bG$ and then neither
for $\bH$, so $P\le\bZ^\circ(\bH)=:\bS'$, a torus. Since $\ell$ divides a
unique $\Phi_d(q)$ this implies that $\bS\le\bS'$, i.e., $\bS$ is the Sylow
$d$-torus of $\bZ^\circ(\bC_\bG^\circ(P))$ and thus uniquely determined by $P$.
Since $P$ is characteristic in $\bS^F$ we have $N_G(P)=N_G(\bS^F)$,
and our claim follows by Proposition~\ref{prop:gen NG(P)}.
\end{proof}

\begin{thm}   \label{thm:picky abelian}
 In the above setting, assume that $\ell$ divides a unique cyclotomic factor
 $\Phi_d(q)$ in the order polynomial of $(\bG,F)$. Then an $\ell$-element
 $x\in G$ is picky if and only if $\bC_G(x)=\bC_G(\bS)$ where $\bS\le\bG$ is a
 Sylow $d$-torus with $x\in\bS^F$.
\end{thm}

\begin{proof}
By Lemma~\ref{lem:cases} we are in the situation of
Proposition~\ref{prop:gen norm syl}. 
Hence $x\in G$ lies in a unique Sylow $\ell$-subgroup of $G$ if and only if it
lies in a unique Sylow $d$-torus $\bS$ of $\bG$, if and only if $\bS$ is the
unique Sylow $d$-torus of $\bC_\bG(x)$, if and only if $\bS\le\bZ(\bC_\bG(x))$,
if and only if $\bC_\bG(x)\le\bC_\bG(\bS)$. As $x\in\bS$ we also have
$\bC_\bG(\bS)\le\bC_\bG(x)$, so $\bC_\bG(\bS)=\bC_\bG(x)$. This implies of
course $\bC_G(\bS)=\bC_G(x)$. Conversely, assume $\bC_G(\bS)=\bC_G(x)$ but
$\bC_\bG(\bS)<\bC_\bG(x)$. Since $\bC_\bG(\bS)$ is a Levi subgroup, hence
connected, this means that either $\bC_\bG(x)$ is disconnected, or of strictly
larger dimension than $\bC_\bG(\bS)$. The first is not possible as $\ell$ is
not a torsion prime. In the second case, as $\bC_\bG(\bS)$ contains a maximal
torus of $\bG$, the maximal unipotent subgroups of $\bC_\bG(x)$ must have larger
dimension than those of $\bC_\bG(\bS)$. But then $|\bC_G(x)|_p>|\bC_G(\bS)|_p$,
a contradiction. In conclusion, $x\in G$ is picky if and
only if $\bC_G(\bS)=\bC_G(x)$, for $\bS$ the (unique) Sylow $d$-torus
containing~$x$.
\end{proof}

\begin{exmp}
 In the setting of Theorem~\ref{thm:picky abelian}, assume $d=1$. Then $x\in G$
 is picky if and only if $x$ is regular. Indeed, in this case the centraliser of
 a Sylow $1$-torus is a maximal torus (namely, a maximally split torus) of
 $\bG$, so the condition becomes: $\bC_G(x)$ is a torus, which means $x$ is
 regular. More generally this characterisation continues to hold whenever
 $d$ is a regular number for $(W,F)$ (in the sense of Springer).
\end{exmp}

The precise determination of subnormalisers in the abelian Sylow case seems
to require a classification of the overgroups of normalisers of Sylow $d$-tori
for simple algebraic groups, a non-trivial and interesting problem, as the
following example shows (see also Corollary~\ref{cor:subn ss} in the algebraic
group case):

\begin{exmp}
 Let $\bG$ be of type $F_4$ and $\ell>3$ dividing $\Phi_3(q)$. Let
 $x\in G=\bG^F$ be an $\ell$-element with centraliser $A_2(q).\Phi_3$. There is
 a subgroup $\tw3D_4(q).3$ of $G$ containing the normalisers of all Sylow
 $d$-tori of $\bC_G(x)$, and hence also $\Sub_G(x)$ by
 Proposition~\ref{prop:gen norm syl}. But it is not immediately obvious that
 the subnormaliser could not be smaller.
\end{exmp}

\subsection{Picky elements in the non-abelian case}

\begin{prop}   \label{prop:Zsigmondy}
 Let $\bS\le\bG$ be a Sylow $e$-torus of $(\bG,F)$ for some $e\ge1$. Then
 the order polynomial of $\bC_\bG(\bS)$ is not divisible by cyclotomic
 polynomials $\Phi_{e\ell^i}$ with $i\ge1$ and $2<\ell$, except when
 $\bG^F=\tw3D_4(q)$, $e\le2$ and $\ell^i=3$.
\end{prop}

\begin{proof}
The centraliser $\bC:=\bC_\bG(\bS)$ is an $F$-stable Levi subgroup of $\bG$.
Its structure is
described in \cite[Exmp.~3.5.14 and~3.5.15]{GM20} for $\bG$ of classical type.
Note that since $\bS$ is a Sylow $e$-torus, it has the form $\bT\bH$ with
$\bT$ a torus whose order polynomial only involves factors $\Phi_d$ with
$d\le 2e$, and $\bH$ a semisimple group of rank less than $e$. The claim can
easily be checked by inspection of the order formulas (e.g.\ in
\cite[Tab.~24.1]{MT11}). For groups of exceptional type, Table~3.3
in \cite{GM20} shows that $\bC$ is itself a torus, with the required property,
unless either $e=4$ and $\bG$ has type $E_7$, but the claim still holds in that
case, or $\bG^F=\tw3D_4(q)$ and $\ell=3$ (where the claim fails).
\end{proof}

\begin{prop}   \label{prop:relWeyl}
 Assume $\bG$ is simple.
 Let $e=e_\ell(q)$ and assume that $\ell$ divides $|W_e|$. If $W_e$ has a
 normal Sylow $\ell$-subgroup then $e\in\{1,2\}$, and either $\ell=2$ and
 $W_e=\fS_2$ or $W(B_2)$, or $\ell=3$ and $W=\fS_3$ or $W(G_2)$.
\end{prop}

\begin{proof}
The relative Weyl groups are described in \cite[Exmp.~3.5.29]{GM20} for $\bG$
of classical type. Namely, they are either symmetric groups $\fS_n$, in which
case $e\in\{1,2\}$, wreath
products $C_d\wr\fS_n$ with $d\in\{e,2e\}$, or certain subgroups of index~2 in
the latter. Note that $\ell>e=e_\ell(q)$ if $\ell\ne2$. The claim follows in
this case from the known normal structure of~$\fS_n$. For $\bG$ of exceptional
type, the occurring relative Weyl groups are given in \cite[Tab.~3]{BMM93}.
No further examples occurs.
\end{proof}

\begin{lem}   \label{lem:relWeyl2}
 Assume $\bG$ is simple and $\Phi_e$ divides the order polynomial of $(\bG,F)$.
 Assume the order of the relative Weyl group $W_e$ of
 a Sylow $e$-torus of $G$ is divisible by $\ell$, where $\ell>2$. Then the
 normaliser of a Sylow $\ell$-subgroup of $W_e$ contains no elements of order
 $\ell r$ for primes $r\equiv1\pmod\ell$.
\end{lem}

\begin{proof}
If $\bG$ is of classical type, then $W_e/\bO_2(W_e)$ is a symmetric group
$\fS_n$ (see \cite[3.5.29]{GM20}), and in $\fS_n$ the Sylow $\ell$-normalisers
are
easily seen not to contain elements of prime order $r>\ell$. For $\bG$ of
exceptional type, it suffices to consider $W(E_8)$ since all relative Weyl
groups are subquotients of this. Now the only element order $\ell r$ for primes
$2<\ell<r$ occurring in $W(E_8)$ is~15, but here $r{\not\equiv}1\pmod\ell$.
\end{proof}

\begin{thm}   \label{thm:non-ab}
 In the above setting assume that $\ell>3$ and that Sylow $\ell$-subgroups of
 $G$ are non-abelian. Then $G$ possesses no picky $\ell$-elements.
\end{thm}

\begin{proof}
Let $P\le G$ be a Sylow $\ell$-subgroup of $G$ and set $e=e_\ell(q)$. By
\cite[Thm~5.16]{Ma07} there is a Sylow $e$-torus~$\bS$ of $(\bG,F)$ such that
the normaliser $\bN_G(\bS)$ contains the normaliser $\bN_G(P)$ of~$P$.
Now $\bL:=\bC_\bG(\bS)$ is an $F$-stable Levi subgroup of~$\bG$. An application
of \cite[Prop.~5.3]{Ma07} to $[\bL,\bL]$ shows that $[\bL,\bL]^F$ is an
$\ell'$-group. Now $\bL/[\bL,\bL]$ is a torus, hence abelian, so as $P$ is
non-abelian by assumption, $\ell$ must divide the order of the relative Weyl
group $W_e=\bN_G(\bS)/\bC_G(\bS)$.

Let $x\in G$ be an $\ell$-element. First assume $x\in\bS_\ell^F$.
By Proposition~\ref{prop:relWeyl}, since $\ell\ge5$,
$W_e$ has more than one Sylow $\ell$-subgroup, so $x$ lies in two different
Sylow $\ell$-subgroups of $\bN_G(\bS)$ and hence of $G$. Now assume $x$ does
not lie in any Sylow $e$-torus, and let $\bT$ be an $F$-stable torus of $\bG$
containing~$x$. Then the order polynomial of $\bT$ must be divisible by a
cyclotomic polynomial $\Phi_{e\ell^i}$ for some $i\ge1$ and hence $|\bC_G(x)|$
is divisible by $\Phi_{e\ell^i}(q)$. Let $r$ be a Zsigmondy primitive prime
divisor of $\Phi_{e\ell^i}(q)$, which exists as $e\ell^i$ is divisible by a
prime at least~5. On the other hand, by Proposition~\ref{prop:Zsigmondy}, the
order polynomial of $\bC_\bG(\bS)$ is not divisible by $\Phi_{e\ell^i}$,
so $|\bC_G(\bS)|$ is prime to~$r$. Thus $W_e$ must contain an element of order
$\ell r$. Now note that $r\equiv1\pmod\ell$. Thus Lemma~\ref{lem:relWeyl2}
together with Lemma~\ref{lem:centr} show that $x$ cannot be picky in
$\bN_G(\bS)$ and so neither in $G$.
\end{proof}

\begin{thm}   \label{thm:non-ab l=3}
 In the above setting assume that Sylow $3$-subgroups of $G$ are non-abelian.
 Then $G$ possesses a picky $3$-element $x$ if and only if one of:
 \begin{enumerate}[\rm(1)]
  \item $G=\SL_3(4),\SU_3(8)$ or $G_2(8)$;
  \item $G=\SU_n(2)$ with $4\le n\le 8$;
  \item $G=\Sp_{2n}(2)$ with $3\le n\le 5$;
  \item $G=\SO_{2n}^+(2)$ with $4\le n\le 5$;
  \item $G=\SO_{2n}^-(2)$ with $4\le n\le 6$; or 
  \item $G=G_2(2)\cong\SU_3(3).2$, $\tw3D_4(2)$, or $F_4(2)$.
 \end{enumerate}
\end{thm}

\begin{proof}
The proof of Theorem~\ref{thm:non-ab} goes through for the prime $\ell=3$
unless either we are in one of the exceptions of \cite[Thm~5.16]{Ma07}, in an
exception of Proposition~\ref{prop:relWeyl}, or if $\Phi_{e3^i}(q)$ has no
primitive prime divisor for some $i\ge1$, where again $e=e_3(q)$. We discuss
these in turn. The exceptions from \cite[Thm~5.16]{Ma07}
are $G=\SL_3(\eps q)$ with $\eps q\equiv 4,7\pmod9$ and $G=G_2(q)$ with
$q\equiv 2,4,5,7\pmod9$. For $G=\SL_3(\eps q)$, by Lemma~\ref{lem:cover} we may
consider $S=\PSL_3(\eps q)$ instead. The Sylow 3-subgroups of $S$ are
elementary abelian of order~9, and the centraliser of an element of order~3 has
structure $((q-\eps)\times(q-\eps)/3).3$. Since the Sylow 3-normalisers have
the form $3^{1+2}.Q_8$ for the appropriate congruences, there cannot exist any
picky 3-elements for $q>4$ by Lemma~\ref{lem:centr}.

Now let $G=G_2(q)$ with $q\equiv 2,4,5,7\pmod9$ and let $\eps\in\{\pm1\}$ with
$q\equiv\eps\pmod3$. The centraliser $C:=\bC_G(t)\cong\SL_3(\eps q)$ of a
3-central element~$t\in G$ contains a Sylow 3-subgroup of $G$. Thus $t$ lies in
several Sylow
3-subgroups of~$C$ (and thus of $G$) unless $q=2$, hence cannot be picky.
There is one further class of non-trivial 3-elements of~$G$, centralised by an
$A_1$-subgroup (see \cite{CR74,EY86}), but the normaliser in~$G$ of a Sylow
3-subgroup (of order~27) does not contain such a subgroup unless again $q=2$.
In the latter case, $G'=\SU_3(3)$ for which all non-trivial 3-elements are
picky by Theorem~\ref{thm:unip picky}, yielding conclusion~(3).
\par
Next, the exceptions of Proposition~\ref{prop:relWeyl} occur precisely for the
groups $\SL_3(\eps q)$, $G_2(q)$ and $\tw3D_4(q)$, since these are the only
groups with a relative Weyl group of a Sylow $e$-torus isomorphic to $\fS_3$ or
$W(G_2)$. If $G=\SL_3(\eps q)$ we have $3|(q-\eps)$ since we assume Sylow
3-subgroups of $G$ are non-abelian, and even $9|(q-\eps)$ by the previous
paragraph. Then, by \cite[Thm~5.16]{Ma07} the normaliser of a Sylow 3-subgroup
$P$ of $G$ lies inside the normaliser $\bN_G(\bS)$ of a Sylow $e$-torus $\bS$.
If $x$ is a 3-element not lying in a conjugate of $\bS_3^F$, then we may
conclude as in the proof of Theorem~\ref{thm:non-ab} that $x$ is not picky.
If $x\in\bS_3^F$ is not regular, its centraliser involves an $A_1$-type group,
hence is not contained in $\bN_G(\bS)=\bS^F.\fS_3$ and again $x$ is not picky
by Lemma~\ref{lem:centr}. Now assume $x\in\bS_3^F$ is regular. If $\bN_G(\bS)$
has more than one Sylow 3-subgroup, again $x$ is not picky. If $\bN_G(\bS)$ has
a unique Sylow 3-subgroup, that is, the Sylow 3-subgroup is normal, then the
elements of order three in $W_e$ need to centralise the $3'$-part of
$\bS^F$ which by inspection forces $\bS^F$ to be a 3-group. Thus $q-\eps$ is a
power of~3 which together with $q\equiv\eps\pmod9$ forces $q=8$, $\eps=-1$.
In the latter case, all regular $x\in \bS_3^F$ are picky by direct computation,
as claimed in (1).

Next, we consider $G=G_2(q)$ with $q\equiv\pm1\pmod9$ (the other congruences
were already discussed above). The argument is now entirely analogous to the
one given for $\SL_3(\eps q)$, and only $G_2(8)$ gives rise to picky elements,
listed in~(1).

Next, if $G=\tw3D_4(q)$ with $q\equiv\eps\pmod3$, $\eps\in\{\pm1\}$, then the
normaliser of a Sylow 3-subgroup is contained in a torus normaliser of the
form $N=(q^3-\eps)(q-\eps).W(G_2)$ (see the discussion in the proof of
\cite[Thm~5.14]{Ma07}). Now elements of order~3 in $W(G_2)$ act like a field
automorphism on the torus of order $q^3-\eps$. In particular, as soon as
$q^2+\eps q+1$ has a primitive prime divisor (necessarily distinct from~3),
a Sylow 3-subgroup is not normal in $N$ and we may conclude as before. This
leaves $q=2$, so $G=\tw3D_4(2)$. Here, by explicit computation, the elements of
order~9 with centraliser $A_1(q).(q^3+1)$ are picky.
\par

Finally, if there exists no primitive prime divisor for $\Phi_{e3^i}(q)$ with
$i\ge1$, then $i=1$, $e=q=2$. If $x\in G$ is a picky 3-element, then it has
order at most~9, as elements of order~27 have centraliser order divisible by
$\Phi_{18}(2)=3^3 19$ in contradiction to Lemma~\ref{lem:relWeyl2}. We now
discuss the various types. For $G$ of classical type, all elements of order~3
lie in a Sylow $2$-torus of $G$ and thus in at least two Sylow 3-subgroups,
arguing as in the proof of Theorem~\ref{thm:non-ab}. So in particular,
$\bC_\bG(x)$ must contain a torus of rational type $\Phi_6$. 
In $\SL_n(q)$ any centraliser order divisible by $q^3+1$ is in fact divisible
by $(q^6-1)/(q-1)$, hence by~7 when $q=2$, contradicting
Lemma~\ref{lem:relWeyl2}. In the remaining groups of classical type, the
primitive 9th roots of unity can occur at most once as an eigenvalue of $x$,
as otherwise the centraliser contains a subgroup $\SU_2(8)$ and hence elements
of order~7. Let $V$ be the submodule of the natural $G$-module spanned by the
eigenspaces for the eigenvalues of~$x$ of order at most~3 and $H$ the induced
classical group on~$V$. Then $H$ has the same type as $G$, respectively type
$\SO^{-\eps}$ if $G$ has type $\SO^\eps$ (since in that case by what we showed
$x$ has exactly all six primitive 9th roots of unity as eigenvalues and the
induced orthogonal group on the corresponding 6-dimensional sum of eigenspaces
is of minus type). Then $\Sub_G(x)$ will contain the subnormaliser of $x|_V$
(of order~3), and by direct computation the latter is all of $H$ if
$H=\SU_3(2),\Sp_6(2),\SO_8^+(2)$ or $\SO_6^-(2)$. As all of these have
order divisible by a prime larger than~3, we see that $x$ cannot be picky
in $\SU_n(2),\Sp_{2n}(2),\SO_{2n}^+(2),\SO_{2n}^-(2)$ when $n\ge9,6,6,7$
respectively. By direct computation, there do exist picky elements of order~9
in all remaining cases.

Similarly, the claim for $G_2(2),\tw3D_4(2),F_4(2)$ and $E_6(2)$ follows by
explicit computation in \GAP. The maximal subgroup $H=Fi_{22}$ of
$S=\tw2E_6(2)$ contains a Sylow 3-subgroup of~$S$. By direct computation, $H$
possesses one class of picky 3-elements, but this is fused in $S$ with one of
the non-picky classes of $H$, so $S$ has no picky 3-element.
The group $G=E_7(2)$ has six classes of elements of order~9, with centralisers
$$\Phi_2^2\Phi_6.\tw2A_2(2)A_1(2),\ \Phi_2.\tw2A_2(8),\ 
  \Phi_2\Phi_6.\tw2A_3(2)A_1(2),$$
$$\Phi_2\Phi_6.\tw2A_4(2),\ \Phi_2\Phi_6.\tw3D_4(2),\ 
  \Phi_2\Phi_6.A_1(8)A_1(2)$$
(private communication of Frank L\"ubeck), while all elements of order~3 lie in
a Sylow 2-torus by \cite{Lue}. Now the normaliser of a Sylow 3-subgroup of~$G$
is contained in the normaliser of a Sylow 2-torus and thus has the form $3^7.N$
where $N$ is the normaliser of a Sylow 3-subgroup of the Weyl group $W(E_7)$,
whose 2-part is $2^3$. Since all of the above centralisers have order divisible
by at least $2^4$, $G$ cannot have picky 3-elements by Lemma~\ref{lem:centr}.
Finally, $G=E_8(2)$ has four classes of elements of order~9, with centralisers
$$\Phi_2.\tw2A_2(8)A_1(2),\ \Phi_2\Phi_6.\tw2A_4(2)A_1(2),\ 
  \Phi_2\Phi_6.\tw2D_5(2),\ \Phi_2\Phi_6.\tw3D_4(2)A_1(2)$$
(again computed by Frank L\"ubeck), and all elements of order~3 have
centralisers of semisimple rank at least~7 by \cite{Lue}. Since the 2-part
of the normaliser of a Sylow 3-subgroup of~$G$ is just $2^4$, we may argue as
before to see that $G$ has no picky 3-elements.
\end{proof}

The case $\ell=2$ is considerably more tricky. It will by all appearance be
even messier than for $\ell=3$, involving Fermat and Mersenne primes, for
example. We will consider it in forthcoming work with M.~Schaeffer Fry.

\subsection{Subnormalisers in Suzuki and Ree groups}

We now assume that $\bG$ is of type $B_2$, $G_2$ or $F_4$ in
characteristic~$p=2,3$, or~2, respectively, and $F$ is a Steinberg endomorphism
of $\bG$ such that $G=\bG^F$ is a Suzuki or Ree group. As before, $\ell$
denotes a prime distinct from the defining characteristic $p$ of $\bG$. We
obtain the analogue of Theorems~\ref{thm:picky abelian} and~\ref{thm:non-ab}:

\begin{thm}   \label{prop:SuzRee}
 In the above setting, let $x\in G$ be an $\ell$-element. Then $x$ is picky if
 and only if $\ell>3$ and $x$ is regular, while $\Sub_G(x)=G$ otherwise.
\end{thm}

\begin{proof}
First assume that $\ell>3$ and $x$ is regular. Then $\ell$ does not divide the
order of the Weyl group of~$\bG$ and hence the Sylow $\ell$-subgroups of~$G$
are abelian (see \cite[Cor.~3.13, p.~259]{BM92}). We may now proceed as in the
proof of Theorem~\ref{thm:picky abelian}, replacing cyclotomic polynomials over
$\QQ$ by suitable cyclotomic polynomials over $\QQ(\sqrt{p})$, with
corresponding Sylow tori, for which the Sylow theorems continue to hold, see
\cite[3.5.3, 3.5.4]{GM20}, to conclude that $x$ is picky.
\par
We now discuss the remaining cases. All semisimple elements $x\ne1$ of
$\tw2B_2(q^2)$ are regular and have order prime to~6, so there is nothing to
prove in this case. In $G={}^2G_2(q^2)$ the only $\ell$-elements $x\ne1$ that
are either non-regular or have $\ell\le3$ are the involutions, with
$\bC_G(x)=\langle x\rangle\times\PSL_2(q^2)$. Such $x$ are also contained in a
Sylow 2-normaliser, of structure $2^3.7.3$, not lying in the maximal subgroup
$\bC_G(x)$. Thus $\Sub_G(x)=G$ as claimed.   \par
Let $G=\tw2F_4(q^2)$. The conjugacy classes of semisimple elements are given
in \cite[Tab.~IV]{Sh75}. First assume $\ell>3$. The non-regular representatives
$x\ne1$ are $t_i$ with $i\in\{1,2,4,5,7,8\}$, where $t_1,t_2,t_5$ and $t_8$
only occur if $q^2\ge8$. The non-cyclic Sylow $\Phi$-tori are $T(j)$ for
$j\in\{1,6,7,8\}$ in \cite[(3.1)]{Sh75}, with relative Weyl groups of order
$16,96,96,48$ respectively. Comparison with the list of maximal subgroups
of~$G$ in \cite[Main Thm]{Ma91} shows that none of them can contain both
$\bC_G(x)$ and $\bN_G(\bT)$ for $\bT$ a Sylow $\Phi$-torus of~$\bG$
containing~$x$.
\par
So finally assume $\ell=3$. Any 3-element $x\in G$ is contained in a Sylow
$\Phi_4$-torus $\bT$ of~$\bG$, where $\bN_G(\bT)=(q^2+1)^2.G_{12}$ with the
primitive complex reflection group $G_{12}$ of order~48. As $\bN_G(\bT)$ is a
maximal subgroup of $G$ by \cite[Main Thm]{Ma91}, either $\Sub_G(x)=\bN_G(\bT)$
or $\Sub_G(x)=G$. 
The structure of a Sylow 3-normaliser $N$ is discussed in \cite[Thm~8.4]{Ma07}.
If $q^2\equiv8\pmod9$ then $N$ is contained in the torus normaliser
$(q^2+1)^2.G_{12}$. Now $G_{12}$ has non-normal Sylow 3-subgroups and hence
arguing as in the proof of Theorem~\ref{thm:non-ab} there cannot exist picky
3-elements, whence $\Sub_G(x)=G$ by what we said above. On the other hand, if
$q^2\equiv2,5\pmod9$ a Sylow 3-normaliser is isomorphic to $\SU_3(2).2$,
not contained in any conjugate of $\bN_G(\bT)$ by \cite[Thm~8.4(b)]{Ma07},
so again we have $\Sub_G(x)=G$.
\end{proof}

By explicit computation, the subnormaliser of an element of order~3 in the Tits
group $\tw2F_4(2)'$ is a maximal subgroup $\PSL_3(3).2$ (not invariant under
the outer automorphism of order~2).

\section{Subnormalisers in algebraic groups}   \label{sec:alg group}
The concept of subnormaliser of course also makes sense in the algebraic
group setting. All algebraic groups considered are over an algebraically closed
field of characteristic $p\ge0$. If $\bG$ is connected reductive, $\bT$ denotes
a maxmial torus, $\Phi$ the root system, $\bU_\al$ for $\al\in\Phi$ the root
subgroups of $\bG$ with respect to $\bT$, and $\bB\ge\bT$ a Borel subgroup
of~$\bG$. Not so surprisingly the situation turns out to be much simpler than
for the finite reductive groups.

\subsection{Subnormalisers of unipotent elements}
We first determine subnormalisers for unipotent elements.

\begin{lem}   \label{lem:uni=nil}
 Let $\bU$ be a unipotent algebraic group and $x\in\bU$. Then
 $\langle x\rangle\subn \bU$.
\end{lem}

\begin{proof}
By \cite[Prop.~2.9]{MT11} we may embed $\bU$ into the (unipotent) group of upper
uni-triangular matrices of $\GL_n$ for a suitable~$n$, so without loss we may
assume $\bU$ is connected. By \cite[Cor.~2.10]{MT11} the unipotent group $\bU$
is nilpotent, so by
\cite[Prop.~17.4]{Hu}, for any proper closed subgroup $\bH<\bU$ we have
$\dim \bH<\dim \bN_\bU(\bH)$. Thus by induction on the dimension, any closed
subgroup is subnormal in $\bU$. As $\langle x\rangle$ is normal in its abelian
(see \cite[Lemma~15.1.C]{Hu}) closure $\overline{\langle x\rangle}$ this
achieves the proof. 
\end{proof}

For $H$ a subgroup of a linear algebraic group we let $\hRu(H)$ denote the
largest normal unipotent subgroup of $H$, a characteristic subgroup.

\begin{lem}   \label{lem:max rad}
 Let $\bG$ be a linear algebraic group, $H\le\bG$ and $x\in\bG$ unipotent.
 If $\langle x\rangle\subn H$ then $x\in\hRu(H)$. In particular, if $H$ is
 maximal (with respect to inclusion) with $\langle x\rangle\subn H$ then
 $H=\bN_\bG(\hRu(H))$.
\end{lem}

\begin{proof}
Let $\langle x\rangle\unlhd H_1\unlhd\cdots\unlhd H_r=H$ be a subnormal
series. Set $U:=\hRu(H)$ and $\tilde H_i:=H_iU$. Since $x$ normalises $U$,
the group $\langle x,U\rangle$ is unipotent. Then
$\langle x,U\rangle\unlhd\tilde H_1$ implies
$\langle x,U\rangle\le\hRu(\tilde H_1)=:U_1$.
As $\tilde H_1\unlhd\tilde H_2$ we now have $U_1\unlhd \tilde H_2$ and thus
$U_1\le\hRu(\tilde H_2)=:U_2$. By induction this yields
$\langle x,U\rangle\le U_1\le\cdots\le U_r=\hRu(\tilde H_r)$,
whence $x\in \hRu(H)$.   \par
Now the Zariski closure $\overline{\hRu(H)}$ is unipotent
\cite[Prop.~2.9]{MT11}, giving $\langle x\rangle\subn\overline{\hRu(H)}$ by
Lemma~\ref{lem:uni=nil} and then of course $\langle x\rangle\subn\hRu(H)$.
Hence, $\langle x\rangle\subn\bN_\bG(\hRu(H))\ge H$ and the last assertion
follows.
\end{proof}

\begin{lem}   \label{lem:unip rad}
 Let $\bG$ be connected reductive. Then the closed unipotent subgroups $\bU$ of
 $\bG$ with $\bU=\hRu(\bN_\bG(\bU))$ are exactly the unipotent radicals of
 the parabolic subgroups of~$\bG$.
\end{lem}

\begin{proof}
Let $\bU\le\bG$ be a closed unipotent subgroup with $\bU=\hRu(\bN_\bG(\bU))$.
By \cite[Cor.~17.15]{MT11}, $\bU$ lies in some Borel subgroup of $\bG$ and then
by \cite[Thm~17.10]{MT11} there is a parabolic subgroup $\bP\le\bG$ with
$\bU\le R_u(\bP)$ and $\bN_\bG(\bU)\le\bP$. In particular, $\bN_\bG(\bU)$
normalises $R_u(\bP)$. Thus, $V:=\bN_{R_u(\bP)}(\bU)$ is contained in and
normalised by $\bN_\bG(\bU)$, so contained in $\hRu(\bN_\bG(\bU))=\bU$, whence
$V=\bU$. By \cite[Prop.~17.4]{Hu} this forces $\bU=R_u(\bP)$.   \par
Conversely, if $\bP\le\bG$ is parabolic, hence connected,
then $\bP/R_u(\bP)$ is connected reductive, and so has no non-trivial normal
unipotent subgroups (since any such would have to be finite, hence central, but
all central elements of connected reductive groups are semisimple, e.g.\ by
\cite[Cor.~8.13(b)]{MT11}), whence $\hRu(\bP)=R_u(\bP)$.
\end{proof}

We obtain the algebraic group analogue of Proposition~\ref{prop:gen NG(P)}:

\begin{prop}   \label{prop:unip alg grp}
 Let $\bG$ be connected and $x\in\bG$ unipotent. Then $\Sub_\bG(x)$ is
 generated by the Borel subgroups of $\bG$ containing~$x$.
\end{prop}

\begin{proof}
By Lemma~\ref{lem:uni=nil}, $\Sub_\bG(x)$ contains all Borel subgroups of $\bG$
containing~$x$. For the converse, assume $\langle x\rangle\subn H$ for some
subgroup $H\le\bG$, where without loss of generality $H$ is maximal with
respect to inclusion. By Lemma~\ref{lem:max rad} then $H=\bN_\bG(\hRu(H))$
and $x\in Q:=\hRu(H)$. Now $Q$ lies in the unipotent radical $\bU$ of some
Borel subgroup of~$\bG$ (\cite[Cor.~17.15]{MT11}). If $Q<\bU$ then
$Q_1:=\bN_\bU(Q)>Q$, and $\langle x\rangle\subn \langle H,Q_1\rangle>H$, a
contradiction, so in fact $Q=\bU$ and $H=\bN_\bG(Q)=\bN_\bG(\bU)$ is closed
(see \cite[Ex.~10.18]{MT11}). Then $\bN_H(\bU)^\circ$ is a Borel
subgroup of $H$, and thus $H^\circ$ is generated by the $H$-conjugates of
$\bN_H(\bU)^\circ$ (see \cite[Thm~6.10]{MT11}), all of which contain $Q$ and
hence $x$. Thus, $H$ is generated by the normalisers (in $H$) of the maximal
unipotent subgroups (of~$H$ and hence of $G$) it contains. The claim follows.
\end{proof}

\begin{thm}   \label{thm:subnorm bG}
 Let $\bG$ be a simple algebraic group in characteristic~$p>0$ and $x\in\bG$ be
 unipotent. Then $\Sub_\bG(x)=\bG$ if and only if $x$ is not regular.
\end{thm}

\begin{proof}
Any regular $x$ is picky by \cite[Prop.~5.1.3]{Ca}. Now let $x\in\bG$ be
non-regular unipotent.
There is nothing to prove if $\bG$ is of type $A_1$. If $\bG$ is of rank~2, the
class representatives for the corresponding finite groups given in
\cite{Sr68,En72,CR74,En76,EY86} show that any non-regular unipotent element has
a $\bG$-conjugate in $\bB$ not involving any simple root element, and thus
$\Sub_\bG(x)=\bG$ by the analogue of Proposition~\ref{prop:unip subnorm}. Now
assume $\bG$ has rank at least~3. Then a suitable conjugate of $x$ in $\bB$ can
be written
as a product of root elements in which at least one simple root $\al$ of $\bG$
does not occur. If $\al$ is an end node, then the observation for rank~2 shows
that $x$ is conjugate to an element in which neither $\al$ nor the simple root
connected to it occurs. Thus, we may assume $\al$ is not an end node. We can
now complete the proof as in Proposition~\ref{prop:GLn subnorm} using
Proposition~\ref{prop:unip alg grp}.
\end{proof}

\subsection{An extension to disconnected groups}
In this subsection we consider a slightly different situation. Let now $\bG$
be the extension of a connected reductive algebraic group~$\bG^\circ$ in
characteristic~$p$ by a graph automorphism $\si$ of order~$p$. We are
interested in unipotent elements in non-trivial cosets of $\bG^\circ$ in $\bG$.
Every coset of $\bG^\circ$ in $\bG$ contains a unique class of \emph{regular
unipotent elements} satisfying analogous properties to the case of connected
groups (see \cite[Prop.~II.10.2]{Sp}).

Let $\bB\le\bG$ be the normaliser in $\bG$ of a $\si$-stable Borel subgroup
$\bB^\circ$ of $\bG^\circ$, with $\si$-stable maximal torus $\bT^\circ$ and
unipotent radical $\bU=\bU^\circ$. Let $W=\bN_{\bG^\circ}(\bT^\circ)/\bT^\circ$
be the Weyl group of $\bG^\circ$; it is normalised by $\si$, and $\si$ acts on
its set of roots $\Phi$ with respect to $\bT^\circ$, permuting its set $\Delta$
of simple roots. We let $\bar\Delta$ denote the set of $\si$-orbits
in~$\Delta$, and we write $\bar S$ for the set of simple reflections of $W^\si$
constructed as in \cite[Lemma~23.3]{MT11}.

\begin{prop}   \label{prop:ext rank 2}
 Let $\bG^\circ$ be of type $A_2,A_3$ or $A_4$ in characteristic~$2$ and $\bG$
 the extension with the non-trivial graph automorphism $\si$ of order~$2$. Let
 $x\in\bG^\circ\si$ be unipotent. Then $\Sub_\bG(x)=\bG$ unless $x$ is
 regular, when $\Sub_\bG(x)=\bC_\bB(\si)$.
\end{prop}

\begin{proof}
Since $\bU$ is a maximal unipotent subgroup in $\bG^\circ$ and normalised
by~$\si$ we may assume $x\in\bU\si$. So $x=x_0\si$ where $x_0$ is a product of
root elements for $\bG^\circ$. If root elements for representatives from all
$\si$-orbits of simple roots occur in $x_0$, then $x$ is regular and hence
contained in a unique Borel subgroup of $\bG$ by \cite[Prop.~II.10.2]{Sp}, so
$\Sub_\bG(x)=\bC_\bB(\si)$. Representatives for the non-regular outer unipotent
classes in types $A_2,A_3$ and $A_4$ are given in \cite[Tab.~4]{Ma93a}. It
transpires that in each case at least one of the given representatives in the
finite group for a class in the algebraic group has the property that its image
under any reflection in $\bar S$ is still in $\bU\si$, and thus $x$ lies in a
unipotent normal subgroup of any~$\bP^\si$, where $\bP$ runs over the
normalisers in $\bG$ of the minimal $\si$-stable standard parabolic subgroups
of~$\bG^\circ$. Since these generate~$\bG$, we have $\Sub_\bG(x)=\bG$.
\end{proof}

\begin{thm}   \label{thm:disc unip picky}
 Let $\bG$ be such that $\bG^\circ$ is simple and of index~$p$ in $\bG$, and
 $x\in\bG^\circ\si$ unipotent. Then $\Sub_\bG(x)=\bG$ unless $x$ is regular,
 when $\Sub_\bG(x)=\bC_\bB(\si)$.
\end{thm}

\begin{proof}
First, by \cite[Prop.~II.10.2]{Sp} the regular unipotent elements
in~$\bG^\circ\si$
are picky, arguing as in the proof of Proposition~\ref{prop:BN}. Now assume
$x\in\bG^\circ\si$ is unipotent non-regular. Again we may assume $x\in\bU\si$.
If $\bC_{\bG^\circ}(\si)$ has rank~1 or~2, we are done by
Proposition~\ref{prop:ext rank 2}, or $\bG^\circ$ is of type $D_4$, $p=3$ and
$\si$ induces triality. In the latter case, inspection of the representatives
given in \cite[Tab.~8]{Ma93a} shows that $x$ can be chosen to lie in the
maximal normal unipotent subgroups of $\bC_\bP(\si)$, where $\bP$ runs over the
normalisers in $\bG$ of the minimal $\si$-stable standard parabolic subgroups of
$\bG^\circ$ and thus $\Sub_\bG(x)=\bG$.   \par
So we may now assume that $p=2$, $\bC_\bG^\circ(\si)$ has rank at least~3, and
that the claim has been shown for all groups of smaller rank.
As $x$ is not regular, by \cite[Prop.~II.10.2]{Sp} it is a product of root
elements not involving representatives from at least one orbit of simple roots
$\al\in\bar\Delta$. Consider the Dynkin diagram of $\bC_{\bG}^\circ(\si)$,
with nodes labelled by $\bar\Delta$. If $\al$ is not an end node, then we can
conclude verbatim by the arguments in the proof of
Proposition~\ref{prop:GLn subnorm}. If $\al$ is an end node, then the image
of~$x$ in the standard Levi subgroup corresponding to $\al$ and the adjacent
node $\beta\in\bar\Delta$ is non-regular, so has a conjugate not involving
a root element for $\beta$ (by \cite[Tab.~4]{Ma93a}) and hence we are reduced
to the previous case since $\bC_\bG^\circ(\si)$ has rank at least~3. Note that
this also holds trivially if this Levi is of type $A_2^2$.
\end{proof}

\subsection{Subnormalisers of semisimple elements}
We return to the setting of connected reductive groups $\bG$ in arbitrary
characteristic and now consider semisimple elements.

What do subnormalisers of semisimple elements look like in connected
reductive groups? They are related to the \emph{$p$-closed subsystems of $\Phi$}
(see \cite[Def.~13.2]{MT11}); as before it is easy and convenient to reduce to
the case of simple groups.

\begin{thm}   \label{thm:subn alg}
 Let $\bG$ be a simple algebraic group and $s\in\bG$ semisimple. Then there
 exists a $p$-closed subsystem $\Psi\subseteq\Phi$ such that 
 $\Sub_\bG(s)=\bN_\bG(\bG(\Psi))$, the normaliser of the subsystem subgroup
 $\bG(\Psi)$ corresponding to $\Psi$, where $\Psi=\emptyset$, $\Psi=\Phi$, or
 $\Psi$ consists of all roots of a fixed length.
\end{thm}

\begin{proof}
Let $s\in\bG$ be semisimple and $\bT$ a maximal torus of $\bG$ containing $s$
(which exists by \cite[Cor.~6.11]{MT11}). Then $\bN_\bG(\bT)$ is contained in
$\Sub_\bG(s)$, so $\bT$ normalises $\Sub_\bG(s)$. Let
$\Psi:=\{\al\in\Phi\mid \bU_\al\le\Sub_\bG(s)\}$. Then, $\Psi$ is a $p$-closed
subsystem of $\Phi$, by definition. Since $\Sub_\bG(s)$ contains $\bN_\bG(\bT)$
and hence representatives for all Weyl group elements, $\Psi$ consists of full
$W$-orbits of roots, that is, of all short, all long, or all roots, or
$\Psi=\emptyset$ (see \cite[Cor.~A.18]{MT11}).  Now $\Sub_\bG(s)^\circ$ is a
subsystem subgroup since it is normalised by the maximal torus $\bT$. According
to \cite[Cor.~13.7]{MT11} this shows that $\Sub_\bG(s)=\bN_\bG(\bG(\Psi))$. 
\end{proof}

\begin{cor}   \label{cor:subn ss}
 In the situation of Theorem~\ref{thm:subn alg}, either $\Sub_\bG(s)=\bG$, or
 one of the following holds:
 \begin{enumerate}[\rm(1)]
  \item $s$ is regular and $\Sub_\bG(s)=\bN_\bG(\bT)$;
  \item $\bG=B_n$ and $\Sub_\bG(s)=D_n$, or $\Sub_\bG(s)=A_1^n$ when $p=2$;
  \item $\bG=C_n$ with $n,p>2$ and $\Sub_\bG(s)=A_1^n$;
  \item $\bG=G_2$ and $\Sub_\bG(s)=A_2$; or
  \item $\bG=F_4$ and $\Sub_\bG(s)=D_4$.
 \end{enumerate}
 Conversely, the cases~(1)--(5) can only occur when $\bC_\bG(s)$ lies in a
 subsystem subgroup of the given type.
\end{cor}

\begin{proof}
By Theorem~\ref{thm:subn alg} we have $\Sub_\bG(s)=\bN_\bG(\bG(\Psi))$ for
some $p$-closed subsystem $\Psi$ of~$\Phi$ consisting of full orbits of roots
under the Weyl group. If $\Psi=\emptyset$ then $\bG(\Psi)=\bT$ and thus
$\Sub_\bG(s)=\bN_\bG(\bT)$. For the rest of the proof assume
$\Psi\ne\emptyset$. We need to understand the possible $\Psi$, where, of
course, we may assume $\Phi$ has two root lengths. The structure of such
subsystems is given in \cite[Tab.~B.2]{MT11}. By
\cite[Thm~13.14 and Prop.~13.15]{MT11} these are $p$-closed precisely under the
conditions listed in the statement. Note that all of these subgroups do indeed
contain the full normaliser of a maximal torus of~$\bG$.   \par
Since $\Sub_\bG(s)$ contains $\bC_\bG(s)$ the additional claim follows.
\end{proof}

\begin{exmp}
 We have not been able to find sufficient conditions for the
 occurrence of cases~(1)--(5) in Corollary~\ref{cor:subn ss}. Let $\bG=\GL_n$
 and $s\in\bG$ a diagonal element with eigenvalues all the $n$th roots of
 unity, where we assume the characteristic of $\bG$ does not divide $n$. Thus
 $s$ is regular semisimple and $\bC_\bG(s)=\bT$, the diagonal maximal torus.
 But a conjugate of $s$ lies in the group of permutation matrices, the Weyl
 group of $\bG$. Hence $\bN_\bG(\bT)$ contains two elements from the $\bG$-class
 of $s$ not conjugate in $\bN_\bG(\bT)$ and thus $\Sub_\bG(s)=\bG$ by
 Theorem~\ref{thm:subn alg}, even though $s$ is regular.
\end{exmp}


\end{document}